# A macro-micro approach to modeling parking

Ziyuan Gu[a,b], Farshid Safarighouzhdi[b], Meead Saberi[b], Taha H. Rashidi[b,*]

[a]*Jiangsu Key Laboratory of Urban ITS, Jiangsu Province Collaborative Innovation Center of Modern Urban Traffic Technologies, School of Transportation, Southeast University, Nanjing 210096, China*
[b]*Research Centre for Integrated Transport Innovation, School of Civil and Environmental Engineering, University of New South Wales (UNSW) Sydney, NSW 2032, Australia*

\* *Corresponding author: rashidi@unsw.edu.au*

**Abstract**

In this paper, we propose a new macro-micro approach to modeling parking. We first develop a microscopic parking simulation model considering both on- and off-street parking with limited capacity. In the microscopic model, a parking search algorithm is proposed to mimic cruising-for-parking based on the principle of proximity, and a parking-related state tracking algorithm is proposed to acquire an event-based simulated data set. Some key aspects of parking modeling are discussed based on the simulated evidence and theoretical analysis. Results suggest (i) although the low cruising speed reduces the network performance, it does not significantly alter the macroscopic or network fundamental diagram (MFD or NFD) unless the cruising vehicles dominate the traffic stream; (ii) distance to park is not uniquely determined by parking occupancy because factors such as cruising speed and parking duration also contribute; and (iii) multiscale parking occupancy-driven intelligent parking guidance can reduce distance to park yielding considerable network efficiency gains. Using the microscopic model, we then extend, calibrate, and validate a macroscopic parking dynamics model with an NFD representation. The demonstrated consistency between the macro- and micro-models permits integration of the two for online parking pricing optimization via model predictive control. Numerical experiments highlight the effectiveness of the proposed approach as well as one caveat. That is, when pricing on-street parking, the road network connected to the alternate off-street parking lots must have sufficient capacity to accommodate the increased parking demand; otherwise, local congestion may arise that violates the homogeneity assumption underlying the macroscopic model.

## 1. Introduction

Vehicles cruising for parking are a disguised source of road congestion incurring substantial externalities. An early study by Axhausen et al. (1994) showed that the parking search time could reach up to 40 percent of the total travel time. More recently, it was showed by Shoup (2006) that the proportion of cruising vehicles in the traffic stream could rise over 70 percent, and by Inci et al. (2017) that the external cruising cost incurred by a car parking for an hour is of the same order of magnitude as the external congestion cost arisen from the trip. The cause is arguably the imbalance between parking demand and supply, i.e., the demand for parking far exceeds the parking availability. To address this imbalance, transportation economists have widely advocated parking fees as a promising instrument because it does not trigger as strong social and political opposition as road tolls[1] (Verhoef et al., 1995). However, to get the price right, understanding the parking dynamics is critical.

While a few mathematical models of parking were proposed (see Section 1.1), limited attempts have been made to explore the validity of these models versus empirical data. This is not surprising as parking data are typically difficult and costly to acquire in reality. Even if they were acquired, such a data set would seldomly be made public. Thus, simulated data naturally become an alternative provided that the simulation is able to mimic the underlying behavior of parking search. In this paper, we first develop such a microscopic parking simulation model, based on which some key aspects of parking modeling are discussed including the effects of low cruising speed, the estimation of distance to park, and the efficiency gains resulting from multiscale parking occupancy-driven intelligent parking guidance. Using the microscopic model, we then extend, calibrate, and validate a macroscopic parking dynamics model with an NFD representation. The demonstrated consistency between the models permits integration of the two resulting in a new macro-micro approach to modeling parking.

---
[1] See Gu et al. (2018a) for an overview of congestion pricing.

*1.1. Literature review*

Parking models and policies have long been a subject of interest especially among transportation economists (Anderson and de Palma, 2004; Arnott et al., 1991; Arnott and Rowse, 1999; Axhausen et al., 1994; Bifulco, 1993; Glazer and Niskanen, 1992; Verhoef et al., 1995). A consensus was somewhat reached that regulatory parking policies, and in particular parking pricing, could improve urban mobility. Underlying some of these studies is one of the earliest mathematical models to describe parking equilibrium, known as the modified bottleneck model inspired by Vickrey (1963). By reformulating the original bottleneck model to consider an additional parking cost component, Arnott et al. (1991) analyzed the spatiotemporal pattern of commuter parking. A key finding was that, even though a location-dependent parking fee is unable to eliminate queueing, it could induce drivers to park in order of decreasing distance from the city center, thereby concentrating arrival times and reducing aggregate schedule delay costs.

The work of Arnott et al. (1991) laid the foundation for using the bottleneck model to analyze commuter parking, which was later extended in various aspects by (i) Zhang et al. (2008) where an integrated morning and evening commute problem was solved considering road tolls and parking fees; (ii) Zhang et al. (2011) where the effects of parking permits and trading were assessed; (iii) Qian et al. (2012) where travelers' parking choices between two parking clusters were considered with different capacity, parking fees, and accessibility to the destination; (iv) Fosgerau and de Palma (2013) where a time-varying parking fee at the workplace was derived; and (v) Liu (2018); Nourinejad and Amirgholy (2018); Zhang et al. (2019b) where a future parking equilibrium problem was solved considering autonomous vehicles. To further account for the physics of hyper-congestion rather than assuming a fixed bottleneck capacity, Liu and Geroliminis (2016) used the macroscopic or network fundamental diagram (MFD or NFD) (see Geroliminis and Daganzo (2008); Mahmassani et al. (2013) for empirical and simulated evidence) to characterize the aggregated network traffic dynamics within the bottleneck model. The authors explicitly considered cruising-for-parking and thus captured its effects on the network performance.

While Vickrey's bottleneck model has clearly attracted much attention, other mathematical tools are also available to analyze the parking equilibrium pattern. One such tool is the classical static traffic assignment model which was extended by many to consider different aspects of parking. Specifically, Zhang et al. (2019a) formulated a joint equilibrium problem of route and parking location choices when commuters travel with autonomous vehicles that can drop them off at the workplace and then automatically head towards a parking lot. The equilibrium problem solved by Boyles et al. (2015) accounted for the mutual dependence between the probabilities of finding parking spots at different locations and the search processes employed by drivers, with the objective of minimizing the expected total travel time or cost. Other factors such as parking duration (Lam et al., 2006), parking capacity (Li et al., 2007), cruising-for-parking flows (Gallo et al., 2011; Leurent and Boujnah, 2014), parking search routes (Pel and Chaniotakis, 2017) were also considered, respectively, in the model formulation. As an equilibrium-seeking model, time-varying cruising-for-parking dynamics is naturally absent as well as its interactions with road congestion.

A few attempts were made to study the interactions between cruising-for-parking and road congestion, mainly from an economic perspective (Arnott and Inci, 2006, 2010; Arnott and Rowse, 2009). The authors proposed a structural model decomposing the parking dynamics system into different pools of vehicles based on their parking-related states. A vehicle can be in transit, cruising for parking, or parked in the system, whose state is updated as soon as it completes the journey in the current pool and joins another. Such a multi-pool flow transferring model can be characterized by a system of mass conservation differential equations, leading to the discovery of both the steady-state solution and the transient dynamics. The same multi-pool representation was recently revisited and extended by Cao and Menendez (2015) using a probabilistic approach. The authors macroscopically modeled the interactions between different pools of vehicles via a parking-related state transition matrix describing how vehicles transition from one state to another over time. Such a probabilistic model was later applied to evaluate different on-street parking policies (Cao et al., 2017; Jakob et al., 2018).

There is a tendency towards modeling on-street parking and/or its interactions with off-street parking (Calthrop and Proost, 2006; Inci and Lindsey, 2015; Najmi et al., 2021; Zheng and Geroliminis, 2016). Limited attempts were made to specifically study the off-street parking pattern or dynamics, possibly due to the fact that off-street parking lots are isolated from the on-road traffic system and thus vehicles traveling therein hardly contribute to road congestion. Nevertheless, a few models were proposed to help optimally operate multiple off-street parking lots. One such model was proposed by Qian and Rajagopal (2014a) on a general parking network consisting of multiple origin-destination pairs connected to different off-street parking lots by driving and walking links. Assuming drivers have perfect knowledge of both the time-varying parking fees and the expected cruising time (by means of either self-experience or online information), the authors found that the system optimum was usually achieved when the more convenient parking cluster was pricing-controlled within an occupancy range from 85 percent to 95 percent. Under the same assumption, the authors also formulated a real-time occupancy-driven parking pricing model as a stochastic control problem (Qian and Rajagopal, 2014b), with the finding that parking pricing and information provision can jointly serve as a dynamic stabilized controller to minimize the total system travel time.

From a system dynamics perspective, capturing the time-varying interactions between cruising-for-parking and road congestion is of importance. Agent-based or microscopic models are expected to excel in this regard because of the vast details these models are able to absorb. One well-known agent-based model was developed by Benenson et al. (2008) using ArcGIS, which was later applied to quantify the effects of cruising-for-parking (Levy et al., 2013; Martens et al., 2010) and various parking policies (Levy et al., 2015; Martens and Benenson, 2008). Other such models include Dieussaert et al. (2009) using a cellular automaton and Waraich and Axhausen (2012) using MATSim. It is, however, intuitive that the more details a model captures, a higher computational cost it usually imposes.

One may thus wonder if parking dynamics could also be reasonably modeled using simplified and computationally efficient approaches. The answer is yes. By using the NFD, a few parsimonious models were recently proposed to macroscopically characterize parking dynamics including cruising-for-parking. The very first model was proposed by Geroliminis (2015) targeting on-street parking only. Analogous to the multi-pool structural model, four families of vehicles, as well as their interactions, were modeled leading to a system of mass conservation differential equations. The NFD was used in conjunction with Little's formula to estimate the time-varying outflow of each family. Such an accumulation-based approach was later extended by (i) Zheng and Geroliminis (2016) to consider a bi-modal network where demand was split by a nested logit model among public transit, parking on and off street; and (ii) Gu et al. (2020) to consider limited on- and off-street parking supplies and interactions between the two. Alternatively, Leclercq et al. (2017) proposed a trip-based approach[2] to macroscopically simulating on-street parking. While characterizing the network-level traffic state by a speed-accumulation NFD, the model individualizes vehicles' travel distances so that user heterogeneity can be incorporated. This is a desirable feature for modeling cruising-for-parking.

---

[2] See Mariotte et al. (2017) for an overview of the accumulation- and trip-based approaches.

*1.2. Objectives and contributions*

It is no surprise that the NFD plays an important role in transportation network modeling. Compared with agent-based or microscopic models, the macroscopic approach requires far less data input and thus simplifies the modeling of complex transportation systems especially at large scale. However, the parsimony of the macroscopic models raises the question on how accurate they are versus micro-simulated or empirical data. This question was recently answered by Paipuri et al. (2019) using microscopic simulation and Mariotte et al. (2020) using real data measurements, but only in networks without parking. When parking dynamics is to be considered, the same question remains open and is arguably more challenging to address given the complexities of cruising-for-parking.

Thus, in this paper, we endeavor to answer two key questions: (i) whether a macroscopic model can accurately capture the detailed parking dynamics resulting from a microscopic model; and (ii) whether the macro- and micro-models can be integrated into a joint approach to modeling parking. Towards this end, we first develop in Section 2 a microscopic parking simulation model that explicitly considers cruising-for-parking. In Section 3, we discuss some key aspects of parking modeling based on the simulated evidence and theoretical analysis. Using the microscopic model, we then extend, calibrate, and validate a macroscopic parking dynamics model with an NFD representation in Section 4. In the same section, we also integrate the macro- and micro-models for online parking pricing optimization. We draw conclusions and provide some further thoughts in Section 5.

Overall, the paper provides the following contribution:

i. A microscopic parking simulation model is developed considering both on- and off-street parking with limited capacity. In the microscopic model, a parking search algorithm is proposed to mimic cruising-for-parking based on the principle of proximity, and a parking-related state tracking algorithm is proposed to acquire an event-based simulated data set.
ii. New insights into various aspects of parking modeling are provided based on the simulated evidence and theoretical analysis, including the effects of low cruising speed, the estimation of distance to park, and the efficiency gains resulting from multiscale parking occupancy-driven intelligent parking guidance.
iii. The microscopic model is used to extend, calibrate, and validate a macroscopic parking dynamics model with an NFD representation. The demonstrated consistency between the models permits integration of the two resulting in a new macro-micro approach to modeling parking.

## 2. A microscopic parking simulation model

Understanding the effects of cruising-for-parking is critical, which can be achieved using either empirical or simulated data. Unfortunately, the former is arguably difficult to acquire and thus was exploited in limited studies via interview- or vehicle-based surveys (Axhausen et al., 1994; Belloche, 2015; Leclercq et al., 2017; Lee et al., 2017; Van Ommeren et al., 2012). Simulated data is a promising alternative provided that the simulation is able to mimic the underlying behavior of parking search.

In this section, we develop such a microscopic parking simulation model considering both on- and off-street parking with limited capacity. The study area is a parking-attracting neighborhood of approximately 0.3 km$^2$ in Sydney, Australia (see Fig. 1a). Both on- and off-street parking are operated by the local government in the area, where most streets have parking supplies and an off-street parking lot is located near the beach. The microscopic model is developed in VISSIM (see Fig. 1b) where the on-street parking spots are treated as real parking spaces and the off-street parking lot is modeled as a virtual container[3]. To account for vehicles' re-departures and heterogeneous parking duration, trip chains are defined assuming individuals' parking duration follows a uniform distribution.

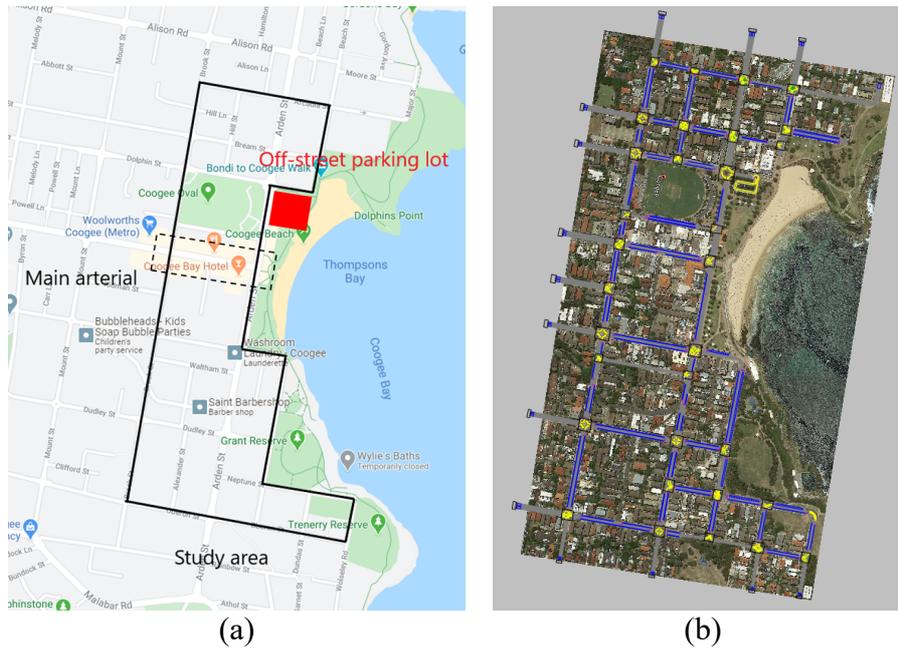

**Fig. 1.** (a) Map of the study area; and (b) bird's eye view of the microscopic model.

In the microscopic model, a parking search algorithm is developed and implemented in real-time (see Fig. 2). Since the neighborhood is treated as a single source of parking attraction[4], we assume drivers are willing to park on any streets with parking supplies and treat them as independent alternatives. This leads to the application of the multinomial logit model[5] for determining the target parking location, whose attributes consist of the parking fee and the parking location attraction. The latter can be thought of as a location-dependent parameter that influences the market share between on- and off-street parking. While the distance or time to park can be implicitly considered by adjusting this parameter, we assume drivers do not have access to such information.

The initial parking location choice may change if that location happens to be fully occupied when the driver arrives. Thus, each link with parking supplies is equipped with a dynamic local parking search decision based on the principle of proximity. The resulting effect is that vehicles would cruise in the vicinity of the target parking locations, referred to as the destinations' search neighborhoods (Benenson et al., 2008; Levy et al., 2015). Specifically, for a typical four-way intersection with parking supplies on all approaches, a driver would randomly choose with equal probability one of the three candidate locations downstream if the target one is fully occupied. Such a decision is made at every intersection along the local search path of the driver until a free spot is found. The "dynamic" feature arises from the fact that the decision is also periodically re-evaluated to reflect the choice heterogeneity of drivers. Note that (i) if none of the downstream links has parking supplies, the candidate location is chosen further downstream at the adjacent intersection; (ii) if the intersection allows U turns, the parking location opposite to the target one is also considered as a candidate; and (iii) if a driver initially intends to park on street, the off-street parking lot is excluded from the candidate set.

---

[3] If the off-street parking lot is not fully occupied, vehicles park without taking up any physical spaces and exit when the parking duration is fulfilled; otherwise, they cruise a circuit to exit the lot and continue the search on street.

[4] Such a single parking destination centroid is a reasonable simplification for a small neighborhood where walking distance does not matter (Cao and Menendez, 2015). But if the size of the neighborhood grows, walking distance may become influential requiring the neighborhood to be partitioned into multiple sub-neighborhoods.

[5] Developing a sophisticated parking location choice model (Hunt and Teply, 1993) is not within the scope of the paper.

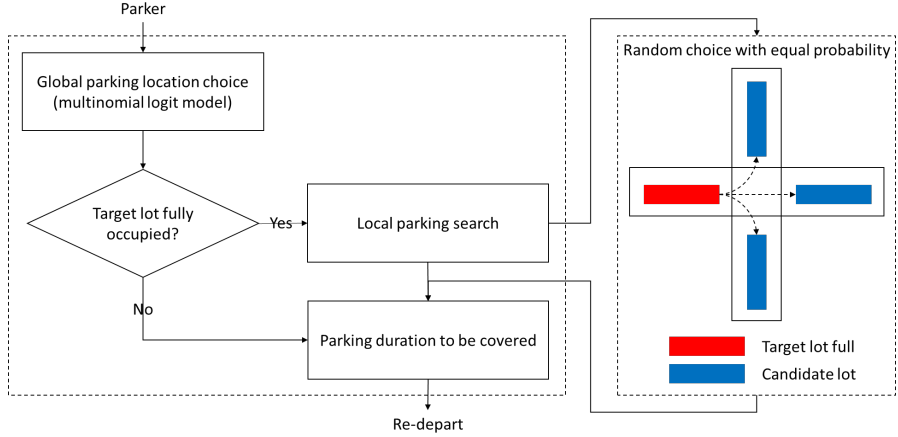

**Fig. 2.** Parking search algorithm in the microscopic model.

One remaining challenge is how different parking-related states of a vehicle can be captured as well as the associated transitions. As this event-based data set is integral to learning the underlying parking dynamics, a parking-related state tracking algorithm is further proposed which communicates with the microscopic model in real-time. The core is to equip each vehicle with a parking-related state vector whose entries are updated whenever a state transition occurs. Given the complexities of the microscopic model, this vector comprises 11 dimensions. Elaboration on each dimension and the detailed algorithmic steps are provided in Section A in the *Supplementary Material*. Implementation of the algorithm requires specifying a time step size which can be set identically to the simulation resolution for achieving the highest numerical accuracy. But the resulting computational time is also the highest. Given this trade-off, we choose to implement the algorithm at every 10 seconds. The obtained event-based data set is then used in the next section to reveal some key aspects of parking modeling.

## 3. Some key aspects of parking modeling

In this section, a few interesting questions on parking modeling are addressed based on the simulated evidence and theoretical analysis.

### 3.1. Effects of low cruising speed

Cruising vehicles cannot be readily recognized from the traffic stream by only judging from the physical appearance. However, from a behavioral perspective, they are admittedly different – they may circulate around certain blocks and travel at much lower speed than the speed limit (Van Ommeren et al., 2012). A question thus arises on the effects of low cruising speed on the network traffic. We emphasize that the extent to which the low cruising speed affects the network traffic depends on the network geometry. If overtaking lanes are available across the network, the effects are arguably minimal as no significant queues would form at these moving bottlenecks. In the subsequent analysis, we only consider links with one traveling lane which is true for most streets in the study area. We first provide some theoretical insights in a simple network and then use the microscopic model to analyze the effects in the study area (which is a more realistic network).

#### 3.1.1. Theoretical analysis in a simple network

Let us consider a square grid network where each link has one traveling lane and one parking lane (see Fig. 3a). If vehicles are homogeneously loaded and not yet cruising for parking, the aggregated network traffic dynamics can be characterized by a speed-accumulation NFD (see the Greenshields model in Fig. 3b). The shape of the NFD would change as soon as one of the vehicles commences cruising at much lower speed, because it creates a local mobile queue that disturbs the homogeneity condition. If the number of cruising vehicles grows such that each link has at least one of them, the network would regain the homogeneity condition as cruising vehicles force the entire traffic stream to move at the same low speed. This is the worst situation possible under the free-flow state. As congestion further builds up in a homogeneous manner, the desired cruising speed cannot be maintained beyond a certain critical density and must decline with the speed of other vehicles towards zero as the network approaches gridlock (see Fig. 3c). The shaded area represents the unstable area of the NFD whose size (i.e. how unstable the network could be) depends on the desired cruising speed. Such network instability only occurs before a certain critical density $k_c = \frac{(v_f - v_c)k_j}{v_f}$, where $k_j$ is the jam density and $v_f$ and $v_c$ are the free-flow and the desired cruising speeds, respectively.

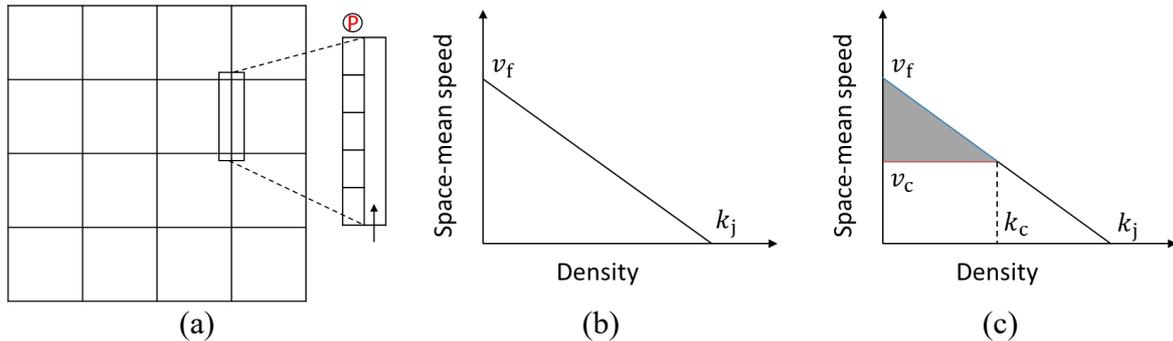

**Fig. 3.** (a) Square grid network with parking supplies; (b) Greenshields speed-density NFD; and (c) hypothetical relationship between density and cruising speed where the shaded area represents the unstable area of the NFD.

What we have yet to answer is how the NFD would change if parking supplies are not evenly distributed in space, and if vehicles are not homogeneously loaded. This is a complex problem involving many possibilities even in the simple square grid network. Thus, to simplify the problem and enable analytical tractability, we use the two-bin idealization (Daganzo et al., 2011) and assume the square grid network is partitioned into two areas where one of them has parking supplies (see Fig. 4a). Each area is modeled as a bin of vehicles using its own NFD (see Fig. 4b). Alternatively, one can think of the two-bin idealization as a network consisting of two links only, where each link is modeled as a bin using the link fundamental diagram.

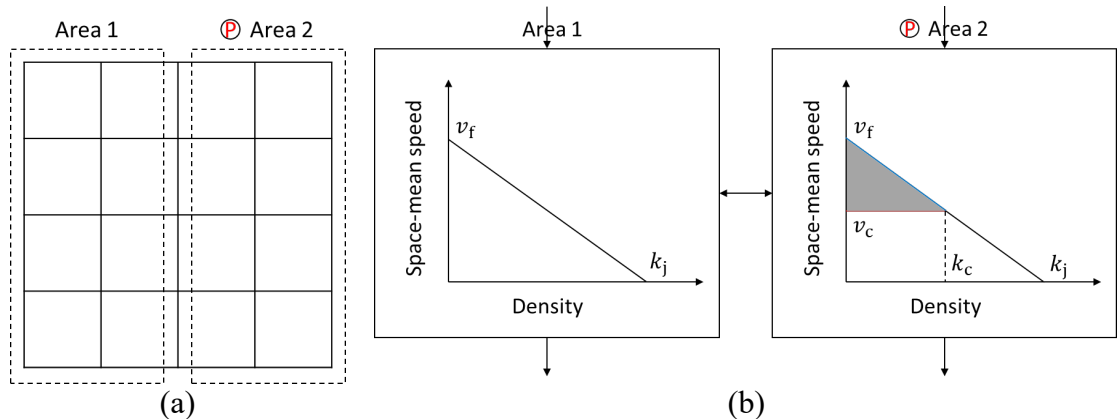

**Fig. 4.** (a) Partitioning the square grid network into two areas; and (b) two-bin idealization.

Let us first derive the NFD of the two-bin system without cruising-for-parking (assuming vehicles entering area 2 are not looking for parking). Let $Q$, $K$, and $V$ denote the space-mean flow, density, and speed of the square grid network, respectively, satisfying the fundamental identity $Q = KV$. Similarly, let $q_x$, $k_x$, and $v_x$ where $x \in \{1,2\}$ denote the space-mean flow, density, and speed of area $x$, respectively, satisfying $q_x = k_x v_x$. Denoting the total network length of each area as $L$,

$$K = \frac{k_1 L + k_2 L}{2L} = \frac{k_1 + k_2}{2}. \tag{1a}$$

Applying Edie's generalized definition of network flow (Edie, 1963; Saberi et al., 2014),

$$Q = \frac{d_1 + d_2}{2LT} = \frac{1}{2}\left(\frac{d_1}{LT} + \frac{d_2}{LT}\right) = \frac{q_1 + q_2}{2}, \tag{1b}$$

where $d_x$ is the total vehicle distance traveled in area $x$ and $T$ is the corresponding time span. Thus,

$$V = \frac{Q}{K} = \frac{q_1 + q_2}{k_1 + k_2} = \frac{k_1 v_1 + k_2 v_2}{k_1 + k_2}, \tag{1c}$$

where $v_1 = v_f - \frac{v_f k_1}{k_j}$ and $v_2 = v_f - \frac{v_f k_2}{k_j}$. Combining Eqs. (1a) and (1c) yields

$$V = \frac{k_1\left(v_f - \frac{v_f k_1}{k_j}\right) + k_2\left(v_f - \frac{v_f k_2}{k_j}\right)}{2K} = v_f - \frac{v_f(k_1^2 + k_2^2)}{2K k_j} = v_f - \frac{v_f}{k_j}\left(2K - \frac{k_1 k_2}{K}\right). \tag{2}$$

While being a function of $K$, $V$ is also affected by how $K$ is distributed between $k_1$ and $k_2$. It can be proved that the upper and lower envelopes have the following forms[6] (see Section B.1 in the *Supplementary Material*):

$$V_{\max} = v_f - \frac{v_f K}{k_j}, \tag{3a}$$

$$V_{\min} = \begin{cases} v_f - \frac{2 v_f K}{k_j}, & K \in \left[0, \frac{k_j}{2}\right] \\ v_f\left(3 - \frac{2K}{k_j} - \frac{k_j}{K}\right), & K \in \left(\frac{k_j}{2}, k_j\right] \end{cases}. \tag{3b}$$

If the two bins are constantly loaded in a homogeneous manner, the NFD of the system is identical to that of each area resulting in the highest possible $V$ for each $K$ (see Fig. 5a). However, if they are loaded as heterogeneously as possible, the system can gridlock at only half of the jam density. When both bins approach the jam density, the upper and lower envelopes naturally converge.

Now, if vehicles entering area 2 are all cruising for parking, we have $v_2 = \min\left\{v_c, v_f - \frac{v_f k_2}{k_j}\right\}$. This is the worst situation for area 2 corresponding to the lower envelope of its own NFD (see Fig. 4b). The upper envelope represents the best situation where no vehicles are cruising and thus is the same as Fig. 5a. We need to consider two conditions in the worst situation, i.e. $k_2 \in [0, k_c]$ and $k_2 \in (k_c, k_j]$, because each condition gives a different $v_2$. It can be proved that synthesizing the upper and lower envelopes derived under each condition yields the following forms (see Section B.2 in the *Supplementary Material*):

---

[6] We are not interested in how cruising and non-cruising vehicles are spatially distributed across the network on a link basis, as there are too many possible combinations. In fact, the NFDs resulting from such combinations would always be bounded by the theoretical upper and lower envelopes.

$$V_{\max} = \begin{cases} v_f - \frac{2v_f K}{k_j}, & K \in \left[0, \frac{k_c}{4}\right] \\ v_c + \frac{(v_f - v_c)k_c}{8K}, & K \in \left(\frac{k_c}{4}, \frac{3k_c}{4}\right] \\ v_c + (v_f - v_c)\left(2 - \frac{k_c}{K}\right)\left(1 - \frac{K}{k_c}\right), & K \in \left(\frac{3k_c}{4}, k_c\right] \\ v_f - \frac{v_f K}{k_j}, & K \in (k_c, k_j] \end{cases} \quad (4a)$$

$$V_{\min} = \begin{cases} v_c, & K \in \left[0, \frac{k_c}{2}\right] \\ v_f - \frac{2v_f K}{k_j}, & K \in \left(\frac{k_c}{2}, \frac{k_j}{2}\right] \\ v_c - \frac{v_c k_j}{2K}, & K \in \left(\frac{k_j}{2}, \frac{k_c + k_j}{2}\right] \\ v_f\left(3 - \frac{2K}{k_j} - \frac{k_j}{K}\right), & K \in \left(\frac{k_c + k_j}{2}, k_j\right] \end{cases} \quad (4b)$$

We consider three typical scenarios based on the relationship between $v_c$ and $\frac{v_f}{2}$. The area shaded by light grey (see Fig. 5b-d) represents the unstable area of the NFD due to the low cruising speed in area 2 as well as the uneven distribution of vehicles between the two areas. Superimposing the upper envelope obtained without cruising-for-parking leads to the areas shaded by dark grey. Results suggest that the size of the unstable area grows if vehicles cruise at lower speed, and that a higher cruising speed reduces the possibility and severity of network instability. As the cruising speed approaches the free-flow speed[7], the unstable area shrinks and converges to Fig. 5a. Results also suggest that the effects of low cruising speed mainly lie in the free-flow and less congested states. When congestion grows homogeneously beyond the critical density, the cruising speed is no different from the speed of other vehicles and thus the NFD exhibits no scatter. Only when vehicles are heterogeneously distributed would scatter appear beyond the critical density, but the effects of low cruising speed are limited.

---

[7] While being theoretically possible, this may not be practical unless with the help of automated parking assistance.

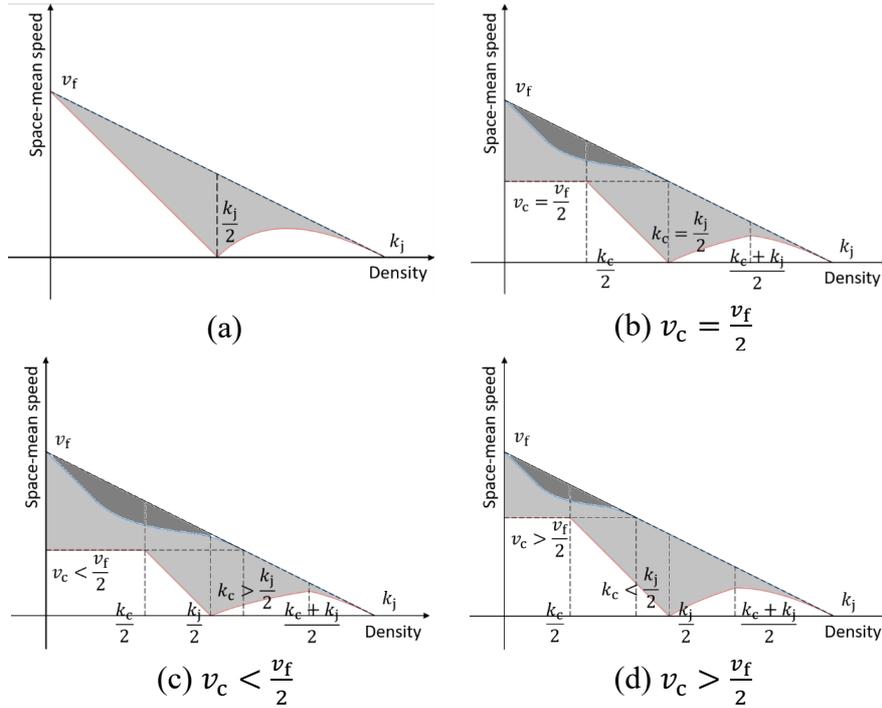

**Fig. 5.** Upper and lower envelopes of the NFD where the shaded area represents the unstable area of the NFD: (a) without cruising-for-parking; and (b)-(d) three typical scenarios with cruising-for-parking.

*3.1.2. Simulated evidence in a realistic network*

Here, we use the microscopic model to analyze the effects of low cruising speed in the study area, where 1,139 parking spots are supplied on street and the off-street parking lot has a capacity of 100. The 1-hr simulation comprises uniformly distributed 1,200 trip chains of drivers seeking to park with different parking duration, which is randomly drawn from a uniform distribution whose average is 30 minutes. While 550 on-street parking spots are constantly occupied throughout the simulation by "captive" users (e.g. long-duration parkers), 360 spots that are already occupied at the beginning of the simulation would be vacated at a constant pace of six spots per minute, representing the re-departures of existing parkers, such as local residents, during the simulation. The simulation also considers vehicles that simply traverse the area without any parking needs (i.e. the passing demand).

In Fig. 6 we compare the simulated NFDs under multiple demand[8] and cruising speed scenarios using Edie's generalized definition of network flow (Edie, 1963; Saberi et al., 2014). Note that drivers' desired speeds are uniformly distributed between $50 \pm 5$ km/hr (the speed limit in the area is 50 km/hr), and that the desired cruising speeds are uniformly and dynamically distributed around a set value (either 10 km/h, 30 km/h, or 50 km/h) also exhibiting $\pm 5$ variations. Results from Fig. 6a-c suggest (i) the NFD enters the congested or gridlocked regime as demand increases; and (ii) the effects of low cruising speed on the NFD appear insignificant. The latter is confusing as it contradicts the theoretical analysis[9]. The cause turns out to be the ratio of the numbers of cruising and non-cruising vehicles. When deriving the theoretical NFD, the two-bin system can be theoretically filled with cruising vehicles only exhibiting the sizable unstable area of the NFD. This is, however, not realistic in the simulation because (i) there exist vehicles traversing the area without any parking needs; and (ii) vehicles seeking to park do not cruise unless they are unable to park at their destinations. Given a predominant proportion of non-cruising vehicles, the effects of low cruising speed on the NFD naturally becomes limited.

---

[8] The passing demand is uniformly distributed within the simulation horizon comprising 1,440, 1,920, and 2,640 vehicles, respectively, corresponding to the low, medium, and high levels.

[9] Although the Greenshields model cannot accurately reflect the "tail" of the simulated NFD, it still reasonably fits the data points when the network is not highly congested. In fact, such highly congested states barely contribute to the theoretical analysis because under such circumstances, the cruising speed is the same as that of the non-cruising traffic. Only when the network is in the free-flow and less congested regimes would we observe speed difference between cruising and non-cruising traffic and thus the effects of the low cruising speed on the NFD.

One may thus wonder if further reducing the parking supplies would alter the NFD. Under such circumstances, we do see more vehicles cruising but the resulting NFD only changes marginally (see the green and blue dots in Fig. 6d). While a more visible change occurs when vehicles are all looking for parking (see the red dots in Fig. 6d), it is still not as notable as the theoretical observations because vehicles do not cruise at low speed unless they fail to park at their destinations. Thus, one can usually be confident of using the simulated or empirical NFD without being concerned about the effects of low cruising speed, except in rare situations where cruising vehicles dominate the traffic stream. As the NFD does not exhibit significant scatter, there is no need for a 3D extension[10] where densities of both cruising and non-cruising vehicles jointly determine the network speed.

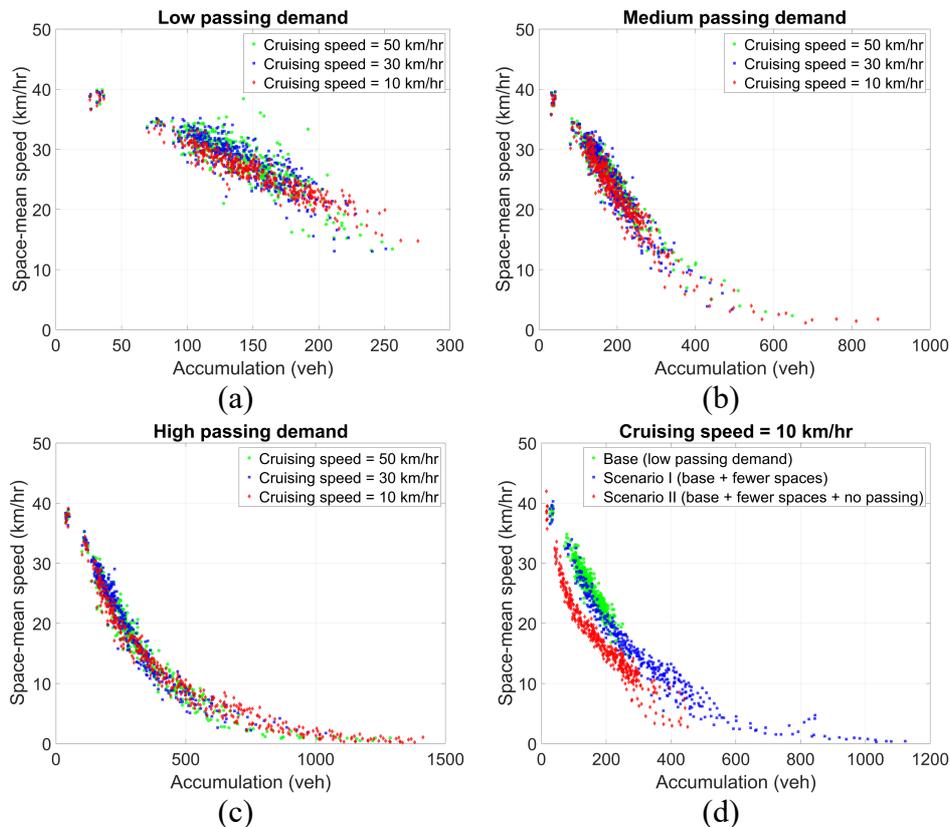

**Fig. 6.** Simulated NFDs under multiple scenarios: (a)-(c) effects of the passing demand with three cruising speeds; and (d) comparing three specific simulation set-ups where the cruising speed is 10 km/hr.

While low cruising speed has limited effects on the NFD, does it change the network performance? We answer this question by comparing four such metrics[11] in Fig. 7. Results suggest that (i) lower cruising speed reduces the network performance resulting in increased average travel time and delay as well as decreased average speed; and (ii) further limiting parking availability deteriorates the network performance. It is counterintuitive, however, that the average travel distance declines as vehicles cruise at lower speed. Given the same parking conditions, one would reasonably envision that, if speed is lower, the distance traveled to find a free spot shall either remain similar under low parking competition or increase under high parking competition. This is true only if vehicles that are already parked do not re-depart. Due to vehicles' re-departures, cruising vehicles are more likely to encounter newly vacated spots at lower speed resulting in shorter distance traveled. This simulated evidence might help explain why in reality vehicles tend to cruise at low speed. Note that when parking availability is limited, drivers might even stop and wait for parked vehicles to re-depart. Under such circumstances, the distance traveled would also decrease (but not the cruising time[12]). This complex behavior is not accounted for in the simulation and thus the above observation is not attributed to it.

---

[10] We refer to Geroliminis et al. (2014) for a 3D-NFD in a bi-modal network consisting of cars and buses.
[11] Due to simulation stochasticity, we consider ten runs with different random seed numbers throughout the paper.

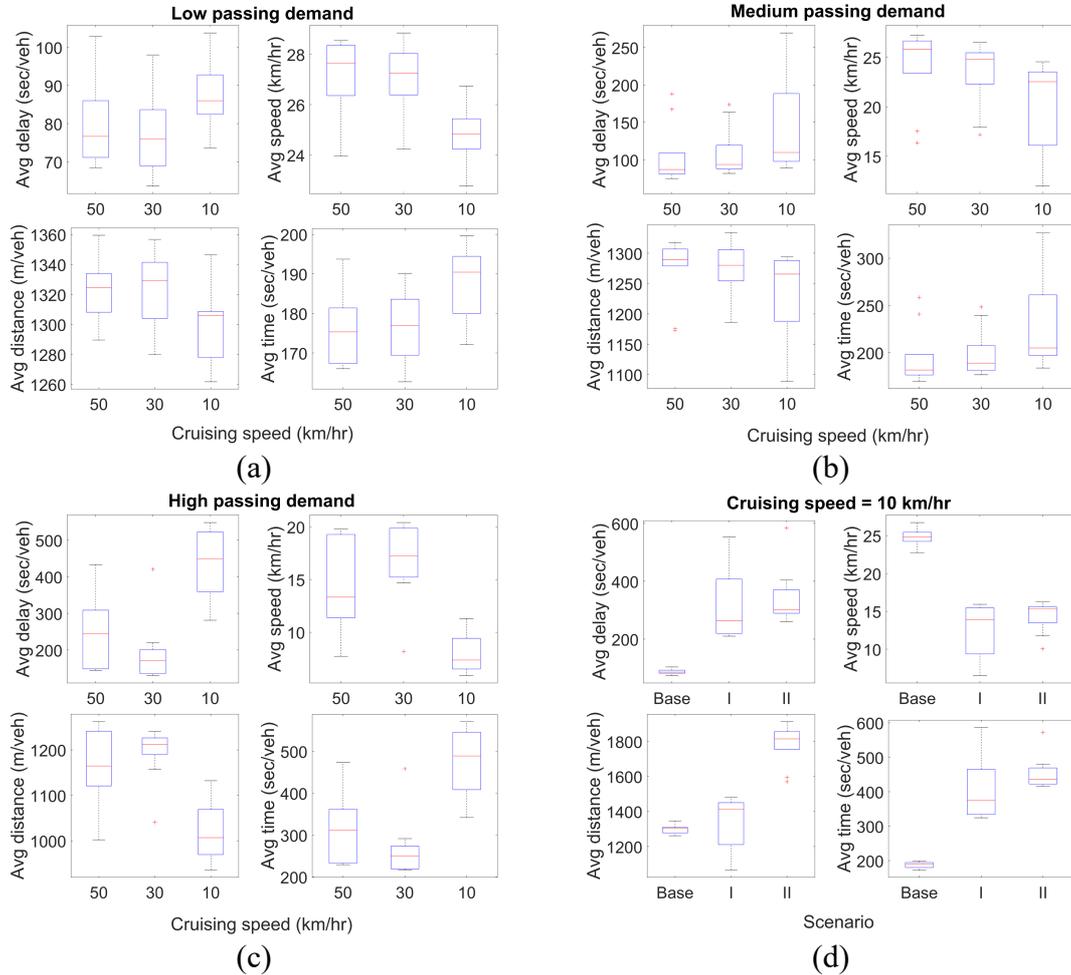

**Fig. 7.** Network performance metrics comparison under multiple scenarios: (a)-(c) effects of the passing demand considering three cruising speeds; and (d) comparing three specific simulation set-ups where the cruising speed is 10 km/hr.

*3.2. Estimation of distance to park*

Distance to park is a key variable in any macroscopic models using the NFD. It is required when applying Little's formula to obtain the outflow of cruising vehicles (see Section 4). In this subsection, we first provide a brief overview of existing methods to estimate distance or time to park and then discuss some new insights gained from the microscopic model.

*3.2.1. A brief overview*

Empirical studies on distance or time to park were often conducted by surveying drivers (Axhausen et al., 1994; Lee et al., 2017; Van Ommeren et al., 2012). Vehicle-based surveys (Belloche, 2015; Leclercq et al., 2017) are a promising alternative that excludes human perception errors, but the implementation is technically and financially challenging. Two functions were typically adopted to characterize the relationship between parking occupancy and time to park:

$$T(t) = a\exp(bO(t)), \tag{5a}$$

$$T(t) = \frac{c}{1-O(t)}, \tag{5b}$$

where $T(t)$ is the time to park that varies with the parking occupancy $O(t)$, and $a$, $b$, and $c$ are parameters to be estimated.

---

[12] Cao and Menendez (2018) recently quantified the cruising time savings through intelligent parking services using the macroscopic model of Cao and Menendez (2015).

Estimating distance to park helps eliminate the mutual dependency between time and speed (Leclercq et al., 2017). It is also estimated as a function of parking occupancy assuming an independent Bernoulli trial for each parking spot screening (Geroliminis, 2015). Let $p(t)$ denote the probability of a spot being free, which is calculated as $p(t) = 1 - O(t)$. Thus, the total number of spots screened $N(t)$ follows a geometric distribution whose probability density function is $P(N(t) = x) = p(t)(1-p(t))^{x-1}$. As the expected value of $N(t)$ is $\frac{1}{p(t)}$, distance to park $L(t)$ can be estimated as

$$L(t) = \frac{d_\text{p}}{p(t)} = \frac{d_\text{p}}{1-O(t)}, \tag{6}$$

where $d_\text{p}$ is the average spacing between two consecutive spots. To ensure the validity of Eq. (6), two requirements must be met: (i) the parking occupancy varies slowly with time; and (ii) parking spots are (approximately) evenly distributed in space. This estimation method, however, was questioned due to its underestimation especially at high parking occupancy (Arnott and Williams, 2017)[13]. To account for the distance traveled on each link where parking is not available (denoted by $d_\text{np}$), Leclercq et al. (2017) proposed and calibrated a modified estimation method:

$$L(t) = \frac{d_\text{np}}{1-O(t)^m} + \frac{d}{1-O(t)}, \tag{7}$$

where $m$ is the average number of parking spots on each link.

*3.2.2. Some new insights*

Distance to park is not a univariate function of parking occupancy. Many other factors are also contributing including parking duration[14] or turnover (Martens et al., 2010). We provide simulated evidence in Fig. 8a-d under four parking duration scenarios. To enable a wider coverage of the simulated data, we increase the number of parkers to 1,800. This is, however, not necessarily the number of observed data points per scenario as some vehicles may still be cruising at the end of the simulation. Results suggest that, when the average parking duration is as low as 15 minutes, the maximal average parking occupancy is only around 80 percent and most distances to park are below one kilometer. However, as soon as the duration doubles, the maximal occupancy approaches 100 percent resulting in as long as four kilometers of distance to park (which quadruples). The situation gets even worse if the duration is extended further. Clearly, when the average parking duration varies, multivaluedness of the average cruising distance occurs even for the same average parking occupancy[15]. That the parking duration affects the distance to park is, in some sense, intuitive – if vehicles park without re-departing (i.e. infinite parking duration), distance to park would rise to infinity in theory.

---

[13] The authors used a computer simulation of vehicles cruising for parking around a circle to provide simulated evidence and discussed several effects contributing to the underestimation that are not accounted for by the geometric distribution (e.g. the bunching effect, the competition effect, and others).
[14] Cruising speed also contributes (see Subsection 3.1.2), but here the focus is on parking duration assuming $v_\text{c} = 30$ km/hr.
[15] Parking duration and parking occupancy are also correlated – a longer parking duration generally results in a faster increase of parking occupancy.

By fitting an exponential function to the mean observations (see Fig. 8e), we see a steeper exponential growth in distance to park with an increasing parking duration. The fitted exponential function is further compared with the two probabilistic estimation methods discussed in the previous subsection (see Fig. 8f). All the fitted curves look similar before the average parking occupancy hits 95 percent, beyond which both the probabilistic methods yield significant overestimation relative to the fitted exponential function. Also, the estimated parameters of the probabilistic methods lose their original physical meanings. Imposing physically meaningful constraints on the parameters only leads to worse goodness of fit (including the aforementioned underestimation especially at high parking occupancy).

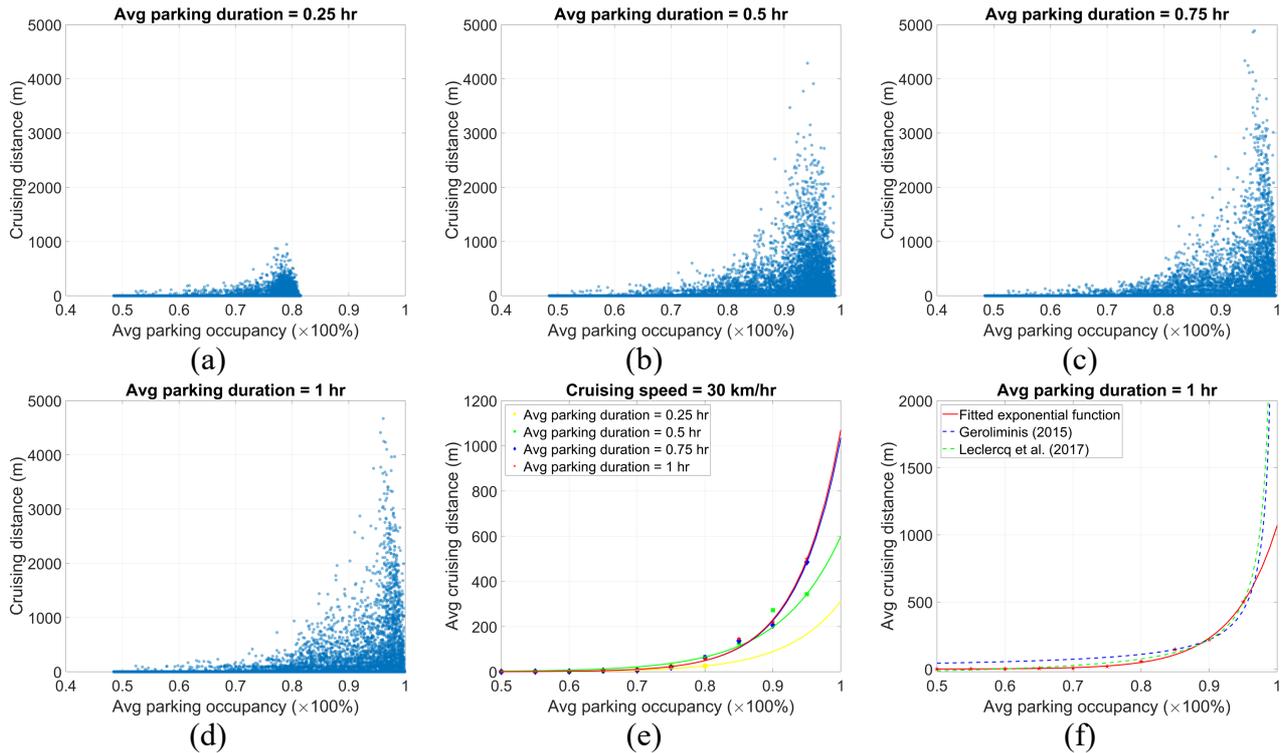

**Fig. 8.** (a)-(d) Simulated observations of distance to park vs. parking occupancy under multiple parking duration scenarios; (e) fitted exponential curves using the mean observations; and (f) comparing three estimation methods.

We also emphasize the effects of the initial parking occupancy that is closely related to the spatial homogeneity/heterogeneity of the parking pattern. At the beginning of the simulation, on-street parking availability of each link, defined as the ratio of the number of free spots and the capacity, is the same so that no link is fully occupied. Thus, as links are all available for parking and vehicles have equal probability of choosing any of them, the on-street parking demand is initially homogeneously distributed. However, once cruising vehicles start to accumulate, the spatial parking pattern grows towards heterogeneity as links with low parking capacity are occupied more quickly giving rise to concentrated cruising traffic in the vicinity. Also contributing to such heterogeneity is the vehicles that fail to park off street. Given the principle of proximity, these vehicles inevitably disturb the homogeneity condition as they cruise nearby the off-street parking lot[16] whenever it is fully occupied.

The spatial heterogeneity of the parking pattern results in longer distance to park. We provide simulated evidence under multiple scenarios with varying initial parking occupancy. Two notable areas are observed in the high parking occupancy range that are mainly filled with the simulated data obtained in the base case (see Fig. 9a). The result suggests that, during the simulation of the parking search from under- to oversaturation, the growing spatial heterogeneity of the parking pattern triggers longer distances to park for certain vehicles. The average distances to park are accordingly lying at the forefront of the occupancy-distance plane (see Fig. 9b). Analogous to the scatter in the NFD due to spatial heterogeneity of congestion distribution (Buisson and Ladier, 2009; Mazloumian et al., 2010)[17], we observe scatter in the relationship between parking occupancy and distance to park due to spatial heterogeneity of the parking pattern. The red upper envelope corresponds to the worst situation where heterogeneity is severest and vehicles cruise the longest on average. In contrast, the black lower envelope represents the best situation where the parking pattern is closest to homogeneity resulting in the least distances to park. Thus, operating on-street parking in a spatially balanced manner is of help in reducing distance to park. This is to be explored further with the help of intelligent parking guidance.

---

[16] We do not consider any tolerance of cruising distance.

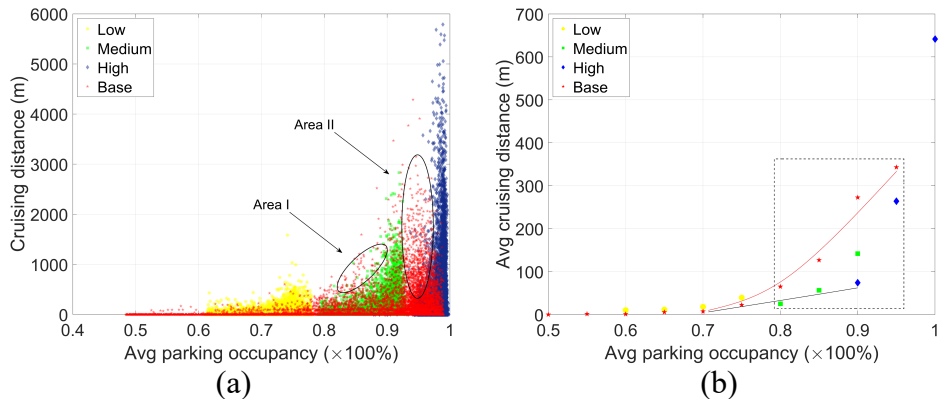

**Fig. 9.** (a) Effects of the initial parking occupancy on distance to park; and (b) mean observations of the simulated data.

*3.3. Efficiency gains from multiscale parking occupancy-driven intelligent parking guidance*

The multiscale parking occupancy-driven intelligent parking guidance consists of an intersection-based local parking search assistance inspired by Leclercq et al. (2017), and a global or regional parking guidance inspired by our previous discussion on Fig. 9. When cruising vehicles approach an intersection without being able to park on the current link, the former assistance provides the real-time parking occupancy information of all the downstream links so that vehicles can be efficiently guided to a free spot. If two or more links have free spots, vehicles would cruise to the link with the maximal probability of parking (i.e. the link with the minimal parking occupancy); if no spot is available on any links, vehicles would randomly cruise to the next intersection for new guidance.

The regional guidance seeks to balance on-street parking operations using only two intersections along the main arterial (see Fig. 1a). The network is thus bi-partitioned where links in the upper region generally have lower parking capacity than those in the lower region. As the off-street parking lot is also located in the upper region, parking competition tends to be severer therein resulting in a heterogeneous parking pattern (e.g. concentration of vehicles that fail to park off street). Thus, the role of the two intersections is to discourage northbound vehicles from entering the upper region whenever the regional parking occupancy exceeds a critical threshold of 95 percent[18]. By diverting some vehicles to the lower region, we hope to regain a close to homogeneous parking pattern so as to reduce distance to park. We emphasize, however, that if no drivers cooperate, such guidance does not have any effects; but if they are cooperative[19], significant efficiency gains would ensue (Gu et al., 2020).

---

[17] See Ramezani et al. (2015) and Gu et al. (2019) for heterogeneity-aware perimeter control and congestion pricing applications, respectively.
[18] See Jakob and Menendez (2019) for assessment of information provision in parking garages only.
[19] Other contributing factors, according to Thompson and Bonsall (1997), include drivers' awareness, level of knowledge, and ability to interpret information, but we assume drivers are aware of and can understand the provided parking guidance.

We compare the simulated NFDs under four scenarios consisting of no-guidance, local guidance, global guidance, and joint local and global guidance (see Fig. 10a). The latter leads to the best network performance as one can readily identify a consistent trend in every metric compared (see Fig. 10b). It is encouraging to see that individual distances to park are generally reduced due to the joint local and global guidance (see Fig. 10c). The positive additive effects are prominent as the resulting maximal distance to park is only about one kilometer, one third of that if no guidance is provided. We also calculate the average distances to park to confirm the above observation (see Fig. 10d). Only under the no-guidance scenario does the average parking occupancy fail to reach 95 percent. In other words, more vehicles are able to park thanks to the provided parking guidance, resulting in a higher parking completion rate (see the bar chart embedded in Fig. 10d). Alternatively, under the no-guidance scenario, more vehicles are still cruising at the end of the simulation given that the total number of vehicles remains the same. We may see even longer distances to park if these vehicles are considered.

Since driver compliance plays an important role in the regional parking guidance, we perform a sensitivity analysis on the percentage of compliant drivers when implementing the joint local and global guidance. It is no surprise that the network performance gets improved as soon as drivers cooperate with the provided regional parking guidance (see Fig. 10e). With only 25 percent of drivers being compliant, the network undergoes a significant decrease in the average travel time, delay, and distance, as well as a significant increase in the average speed. Increasing the percentage further does not lead to an even more notable improvement in the network performance, because the threshold to activate the regional parking guidance remains fixed at 95 percent without accounting for the time-varying nature of the critical parking occupancy (Qian and Rajagopal, 2014b). Nevertheless, with more and more cooperative drivers, the average distance to park becomes shorter and shorter (see Fig. 10f).

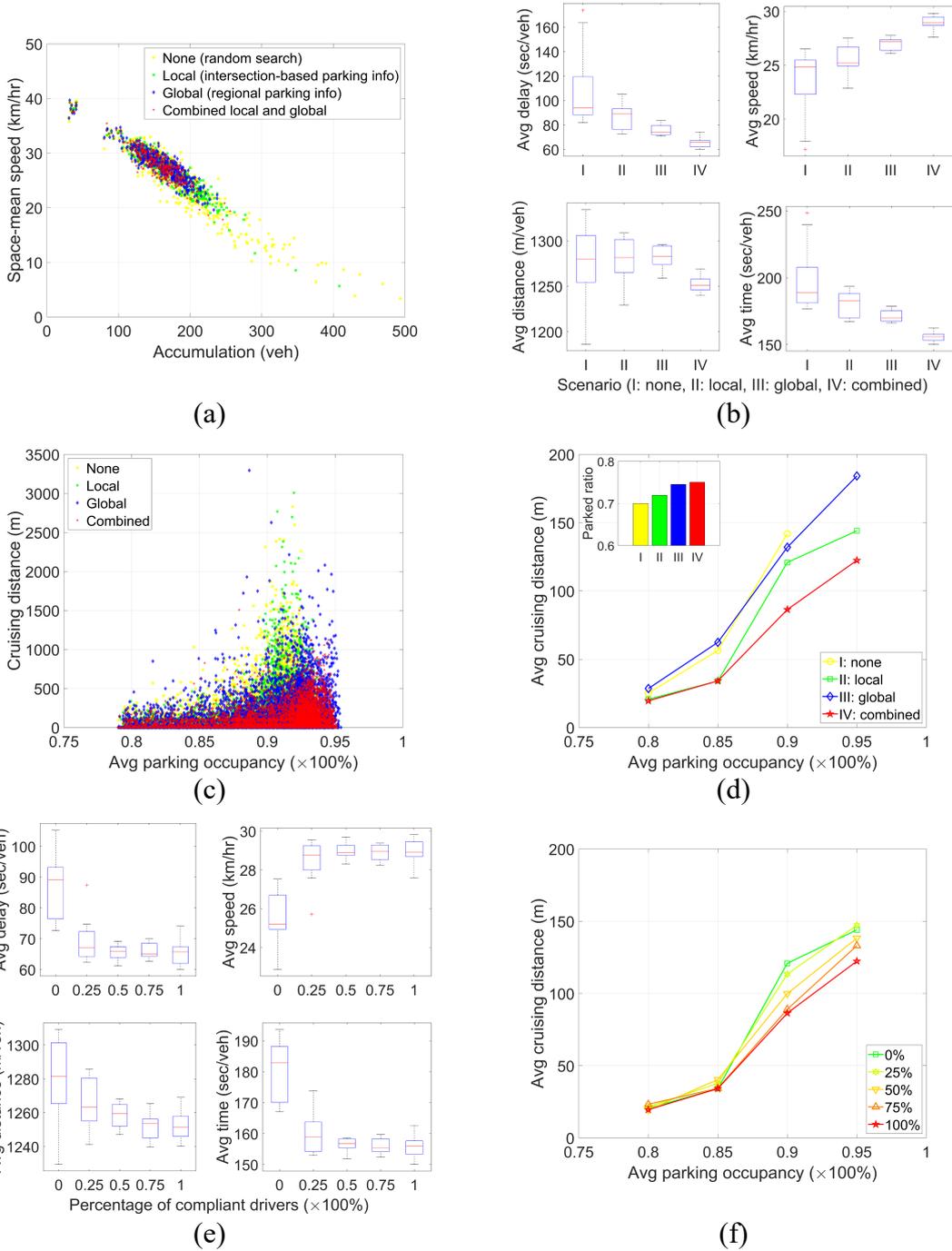

**Fig. 10.** Effects of intelligent parking guidance on the network performance: (a) simulated NFDs; (b) network performance metrics comparison; (c) individual and (d) average distances to park as a function of the average parking occupancy; and (e) and (f) sensitivity analysis on the percentage of compliant drivers with the regional parking guidance.

## 4. A macroscopic parking dynamics model with an NFD representation

Using the microscopic model, we extend, calibrate, and validate a macroscopic parking dynamics model with an NFD representation, and then integrate the macro- and micro-models for online parking pricing optimization.

*4.1. Model development*

The macroscopic model adopts an accumulation-based approach. To ensure model consistency for calibration and validation, it extends previous studies (Geroliminis, 2015; Gu et al., 2020; Zheng and Geroliminis, 2016) in one or more of the following aspects: (i) both on- and off-street parking with limited capacity are considered as well as their interactions; (ii) time delays associated with vehicles transferring from off- to on-street parking search are accounted for; and (iii) a distribution of parking duration is incorporated to model vehicles' re-departures.

In Fig. 11 we illustrate the structure of the multi-pool representation of the model describing flow transfers between different families of vehicles (see Table 1 for nomenclature). Six families co-exist with their unique parking-related states. Specifically, family i consists of vehicles still traveling towards their target parking lots and not yet cruising on street. When these vehicles fail to park at the destinations, they commence the local parking search and join family iv. Vehicles never stay permanently in this family because, as soon as a free spot is found, they park and join family v. Similarly, family ii comprises vehicles still traveling towards the off-street parking lot. Upon reaching the destination, these vehicles park and join family vi if the lot is not fully occupied; otherwise, they cruise out of the lot to continue the search on street and join family iv. In-transit vehicles are part of family iii. The other part are the re-departed vehicles from both on- and off-street parking.

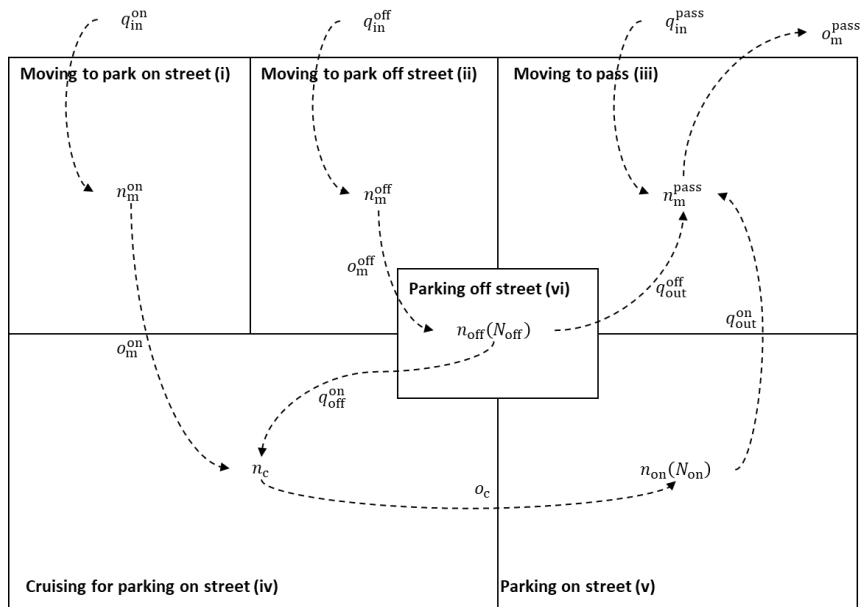

**Fig. 11.** Multi-pool representation of the model with flow transfers between different families of vehicles.

**Table 1**
Summary of variables and parameters.

| Variable | Unit | Description |
|---|---|---|
| $k$ | | Time step between $t_{k-1}$ and $t_k$ |
| $q_{\text{in}}^x(k), x \in \{\text{on, off, pass}\}$ | veh/$\Delta t$ | Vehicles' arrivals to park on or off street, or to traverse the area without any parking needs |
| $q_{\text{out}}^x(k), x \in \{\text{on, off}\}$ | veh/$\Delta t$ | Vehicles' re-departures from on or off street |
| $q_{\text{off}}^{\text{on}}(k)$ | veh/$\Delta t$ | Vehicles that continue the search on street due to a fully occupied off-street parking lot |
| $n_{\text{m}}^x(t_k), x \in \{\text{on, off, pass}\}$ | veh | Accumulation of vehicles that are still moving, either to park on or off street or to traverse the area without any parking needs |
| $n_x(t_k), x \in \{\text{m, c, on, off}\}$ | veh | Accumulation of vehicles that are moving, cruising, or parked on or off street |
| $n(t_k)$ | veh | Accumulation of active (i.e. non-parked) vehicles |
| $P_x(t_k), x \in \{\text{m, c}\}$ | veh·km/hr | Production of moving or cruising vehicles |
| $o_{\text{m}}^x(k), x \in \{\text{on, off, pass}\}$ | veh/$\Delta t$ | Outflow of moving vehicles that commence cruising on street, enter the off-street parking lot, or exit the network |
| $o_c(k)$ | veh/$\Delta t$ | Outflow of cruising vehicles that park on street |
| $O_{\text{on}}(t_k)$ | % | Average on-street parking occupancy |
| $l_c^{\text{on}}(O_{\text{on}}(t_k))$ | km | Expected distance to park on street |
| $v(n(t_k))$ | km/hr | Space-mean speed of the network |
| $v_{\text{on}}(n(t_k))$ | km/hr | Space-mean speed of cruising |
| $\tau_x(k), x \in \{\text{on, off}\}$ | $ | On- or off-street parking price |
| **Parameter** | | |
| $T$ | hr | Simulation horizon |
| $\Delta t$ | hr | Time step size |
| $N_x, x \in \{\text{on, off}\}$ | veh | On- or off-street parking capacity |
| $l_{\text{m}}^x, x \in \{\text{on, off, pass}\}$ | km | Average moving distance to on- or off-street parking, or to traverse the area without any parking needs |
| $v_x^f, x \in \{\text{on, off}\}$ | km/hr | Desired cruising speed on or off street |
| $l_{\text{off}}$ | km | One circuit of the off-street parking lot |
| $k_{\text{off}}$ | | Number of time steps to cruise out of the off-street parking lot |

### 4.1.1. Arrivals and re-departures

Vehicles' arrivals represent the exogenous demand consisting of $q_{\text{in}}^{\text{on}}$, $q_{\text{in}}^{\text{off}}$, and $q_{\text{in}}^{\text{pass}}$. The former two are determined by the choice model while the latter corresponds to the medium demand described in Section 2. We do not account for (i) demand elasticity whereby vehicles' arrivals may be affected by the changing parking conditions; and (ii) boundary inflow or outflow restrictions which may result in local queues (see Mariotte and Leclercq (2019) for effective inflows and outflows).

Vehicles' re-departures represent the endogenous demand when their parking duration is fulfilled. Once these vehicles re-depart, they join family iii to exit the network. To model vehicles' re-departures, we employ a probabilistic approach where a distribution of parking duration is explicitly considered (Cao and Menendez, 2015), as opposed to identical parking duration (Zheng and Geroliminis, 2016) or an exogenously given re-departure rate (Geroliminis, 2015; Gu et al., 2020; Leclercq et al., 2017). Let $f(x)$ denote the probability density function of the distribution of parking duration. To calculate $q_{\text{out}}^{\text{on}}$ and $q_{\text{out}}^{\text{off}}$ in the macroscopic model for any time steps, we use the fact that, as soon as a vehicle is parked, the probability of re-departing during each subsequent time step is an integral of $f(x)$ over all the possible parking duration falling into this time step. In other words,

$$q_{\text{out}}^{\text{on}}(k) = \sum_{i=2}^{k} o_c(i-1) \int_{(k-i)\Delta t}^{(k-i+1)\Delta t} f(x)dx, \tag{8a}$$

$$q_{\text{out}}^{\text{off}}(k) = \sum_{i=2}^{k} \left( o_{\text{m}}^{\text{off}}(i-1) - q_{\text{off}}^{\text{on}}(i-1) \right) \int_{(k-i)\Delta t}^{(k-i+1)\Delta t} f(x)dx. \tag{8b}$$

If the cumulative distribution function $F(x)$ is available, the integral $\int_{(k-i)\Delta t}^{(k-i+1)\Delta t} f(x)dx$ is simply $F\big((k-i+1)\Delta t\big) - F\big((k-i)\Delta t\big)$, which calculates the probability that parking duration is within the range between $(k-i)\Delta t$ and $(k-i+1)\Delta t$. In Eq. (8a), this probability is multiplied by the number of vehicles that successfully park on street during the time step $i-1$. Summing up this product up to the $(k-1)^{\text{th}}$ time step yields the expected number of vehicles that re-depart during the $k^{\text{th}}$ time step. While the same logic applies to Eq. (8b), we must be aware that, when the off-street parking lot is fully occupied, vehicles are unable to park off street and thus continue the search on street. These vehicles $q_{\text{off}}^{\text{on}}(i-1)$ are subtracted from the arrivals $o_{\text{m}}^{\text{off}}(i-1)$ to obtain the true number of vehicles that successfully park off street.

In consistency with the microscopic model, we assume $f(x)$ corresponds to a uniform distribution between zero and one hour. This assumption, however, does not necessarily suggest that the distribution of parking duration is best described by a uniform distribution. From a practical perspective, empirical evidence can be of help in determining the appropriate distribution, which may vary from place to place and from time to time. Nevertheless, given the genericity of the probabilistic approach, other distributions (Cao and Menendez, 2013; Richardson, 1974) can also be used in Eqs. (8a) and (8b). Since the uniform distribution suggests that the probability of re-departing during each time step becomes identical equating to $\frac{\Delta t}{T}$[20], the two equations can be further simplified:

$$q_{\text{out}}^{\text{on}}(k) = \frac{\Delta t}{T} \sum_{i=2}^{k} o_{\text{c}}(i-1), \tag{9a}$$

$$q_{\text{out}}^{\text{off}}(k) = \frac{\Delta t}{T} \sum_{i=2}^{k} o_{\text{m}}^{\text{off}}(i-1) - q_{\text{off}}^{\text{on}}(i-1). \tag{9b}$$

*4.1.2. System dynamics*

Given family-specific accumulations, the total number of active (i.e. non-parked) vehicles is

$$n(t_k) = n_{\text{m}}^{\text{off}}(t_k) + n_{\text{m}}^{\text{on}}(t_k) + n_{\text{m}}^{\text{pass}}(t_k) + n_{\text{c}}(t_k), \tag{10}$$

which is used to calculate the space-mean speed $v(n(t_k))$. The average on-street parking occupancy is calculated as

$$O_{\text{on}}(t_k) = \frac{n_{\text{on}}(t_k)}{N_{\text{on}}}, \tag{11}$$

which is used to calculate the expected distance to park using the exponential function:

$$l_{\text{c}}^{\text{on}}\big(O_{\text{on}}(t_k)\big) = a\exp\big(bO_{\text{on}}(t_k)\big), \tag{12}$$

where $a$ and $b$ are to be estimated.

To calculate network production (Edie, 1963), we distinguish between moving and cruising vehicles due to their difference in speed:

$$P_{\text{c}}(t_k) = n_{\text{c}}(t_k)v_{\text{on}}\big(n(t_k)\big), \tag{13a}$$

---

[20] The simulation horizon $T$ is one hour, thus being equal to the maximal parking duration.

$$P_{\mathrm{m}}(t_k) = n(t_k)v\big(n(t_k)\big) - P_{\mathrm{c}}(t_k), \tag{13b}$$

where $v_{\mathrm{on}}(t_k) = \min\{v_{\mathrm{on}}^{\mathrm{f}}, v(t_k)\}$ and $v(t_k)$ is to be estimated. The family-specific outflows are thus calculated by means of Little's formula (Little, 1961):

$$o_{\mathrm{c}}(k) = \frac{P_{\mathrm{c}}(t_{k-1})\Delta t}{l_{\mathrm{c}}^{\mathrm{on}}(O_{\mathrm{on}}(t_{k-1}))}, \tag{14a}$$

$$o_{\mathrm{m}}^{x}(k) = \frac{P_{\mathrm{m}}(t_{k-1})n_{\mathrm{m}}^{x}(t_{k-1})\Delta t}{l_{\mathrm{m}}^{x}\sum_x n_{\mathrm{m}}^{x}(t_{k-1})}, \quad x \in \{\mathrm{on}, \mathrm{off}, \mathrm{pass}\}, \tag{14b}$$

where $l_{\mathrm{m}}^{x}$ is to be estimated. If the demand structure were to change once for the entire time horizon, a re-estimation of $l_{\mathrm{m}}^{x}$ is required (Batista et al., 2019). If the demand structure were to change over time, then re-estimation could be applied to each time interval where the structure remains unchanged.

The overall system dynamics can be macroscopically characterized by the following system of mass conservation equations in the discrete form (see Fig. 11 for graphical interpretation):

$$n_{\mathrm{m}}^{\mathrm{off}}(t_k) = n_{\mathrm{m}}^{\mathrm{off}}(t_{k-1}) + q_{\mathrm{in}}^{\mathrm{off}}(k) - o_{\mathrm{m}}^{\mathrm{off}}(k), \tag{15a}$$

$$n_{\mathrm{m}}^{\mathrm{on}}(t_k) = n_{\mathrm{m}}^{\mathrm{on}}(t_{k-1}) + q_{\mathrm{in}}^{\mathrm{on}}(k) - o_{\mathrm{m}}^{\mathrm{on}}(k), \tag{15b}$$

$$n_{\mathrm{m}}^{\mathrm{pass}}(t_k) = n_{\mathrm{m}}^{\mathrm{pass}}(t_{k-1}) + q_{\mathrm{in}}^{\mathrm{pass}}(k) + \sum_{x \in \{\mathrm{on,off}\}} q_{\mathrm{out}}^{x}(k) - o_{\mathrm{m}}^{\mathrm{pass}}(k), \tag{15c}$$

$$n_{\mathrm{c}}(t_k) = n_{\mathrm{c}}(t_{k-1}) + q_{\mathrm{off}}^{\mathrm{on}}(k - k_{\mathrm{off}}) + o_{\mathrm{m}}^{\mathrm{on}}(k) - o_{\mathrm{c}}(k), \tag{15d}$$

$$n_{\mathrm{off}}(t_k) = n_{\mathrm{off}}(t_{k-1}) + o_{\mathrm{m}}^{\mathrm{off}}(k) - q_{\mathrm{off}}^{\mathrm{on}}(k) - q_{\mathrm{out}}^{\mathrm{off}}(k), \tag{15e}$$

$$n_{\mathrm{on}}(t_k) = n_{\mathrm{on}}(t_{k-1}) + o_{\mathrm{c}}(k) - q_{\mathrm{out}}^{\mathrm{on}}(k), \tag{15f}$$

where

$$k_{\mathrm{off}} = \left[\frac{l_{\mathrm{off}}}{v_{\mathrm{off}}^{\mathrm{f}}\Delta t}\right]^{21}, \tag{16}$$

$$q_{\mathrm{off}}^{\mathrm{on}}(k) = \max\left\{0, \min\{o_{\mathrm{m}}^{\mathrm{off}}(k), n_{\mathrm{m}}^{\mathrm{off}}(t_{k-1}) + q_{\mathrm{in}}^{\mathrm{off}}(k)\} - \left(N_{\mathrm{off}} - n_{\mathrm{off}}(t_{k-1}) + q_{\mathrm{out}}^{\mathrm{off}}(k)\right)\right\}. \tag{17}$$

Eq. (16) calculates the number of time steps required to cruise out of the off-street parking lot if it is fully occupied. It represents the "time delay" associated with vehicles that reappear on street to continue the search. The number of these vehicles is calculated by Eq. (17) as the difference between two quantities $\min\{o_{\mathrm{m}}^{\mathrm{off}}(k), n_{\mathrm{m}}^{\mathrm{off}}(t_{k-1}) + q_{\mathrm{in}}^{\mathrm{off}}(k)\}$ and $N_{\mathrm{off}} - n_{\mathrm{off}}(t_{k-1}) + q_{\mathrm{out}}^{\mathrm{off}}(k)$. The former calculates the number of vehicles entering the off-street parking lot which is upper bounded by the maximal family-specific accumulation, while the latter calculates the number of free spots in the off-street parking lot considering those vacated by newly re-departed vehicles. The difference is lower bounded by zero when all the vehicles entering the off-street parking lot are able to park.

Lastly, the outflow of each family shall be physically upper bounded by the family-specific accumulation or capacity. Thus, the following inequalities must hold:

---

[21] $[x]$ means the nearest integer to $x$.

$$o_{\text{m}}^{\text{off}}(k) \leq n_{\text{m}}^{\text{off}}(t_{k-1}) + q_{\text{in}}^{\text{off}}(k), \tag{18a}$$

$$o_{\text{m}}^{\text{on}}(k) \leq n_{\text{m}}^{\text{on}}(t_{k-1}) + q_{\text{in}}^{\text{on}}(k), \tag{18b}$$

$$o_{\text{m}}^{\text{pass}}(k) \leq n_{\text{m}}^{\text{pass}}(t_{k-1}) + q_{\text{in}}^{\text{pass}}(k) + \sum_{x \in \{\text{on,off}\}} q_{\text{out}}^{x}(k), \tag{18c}$$

$$o_{\text{c}}(k) \leq \min\{n_{\text{c}}(t_{k-1}) + q_{\text{off}}^{\text{on}}(k - k_{\text{off}}) + o_{\text{m}}^{\text{on}}(k), N_{\text{on}} - n_{\text{c}}(t_{k-1}) + q_{\text{out}}^{\text{on}}(k)\}. \tag{18d}$$

*4.2. Model calibration and validation*

*4.2.1. Calibration*

Compared with the fine-grained microscopic model, the macroscopic model is parsimonious in nature involving several simplifications. A question thus arises on the capability of the latter in reproducing similar results to those obtained from the former. To answer this question, we calibrate and validate the macroscopic model using the microscopic model as the baseline. Note that the microscopic model itself is not calibrated or validated (a limitation to be addressed if parking data become available in the future), but this does not affect the validity of the model's cruising-for-parking simulation and the consistency between the two models to be demonstrated later.

Calibrating the macroscopic model requires that we estimate the following three functions or parameters: (i) speed-accumulation NFD $v(n(t_k))$; (ii) family-specific average moving distances $l_{\text{m}}^{x}, x \in \{\text{on, off, pass}\}$; and (iii) the relationship between parking occupancy and distance to park $l_{\text{c}}^{\text{on}}(O_{\text{on}}(t_k))$. We estimate the simulated NFD using the weighted least square method that addresses the sample selection bias problem (Qu et al., 2015). While different functional forms of the NFD exist, the three-parameter logistic model (Wang et al., 2011) provides the best fit (see Fig. 12a) which reads $v(n(t_k)) = \frac{55.2}{1+\exp\left(\frac{n(t_k)-151.2}{142.1}\right)}$ km/hr.

To estimate family-specific average moving distances, we categorize vehicles' travel distances obtained in the microscopic model according to their memberships in families i-iii. If a vehicle changes the parking-related state along the trip, its entire travel distance is segmented and shared by more than one family. As illustrated in Fig. 12b, the average moving distances for different families do not remain the same across multiple replications due to the microsimulation stochasticity, but the associated variations are rather low. After averaging across all the replications, we end up with $l_{\text{m}}^{\text{off}} = 0.9$ km, $l_{\text{m}}^{\text{on}} = 1.0$ km, and $l_{\text{m}}^{\text{pass}} = 1.1$ km.

The relationship between parking occupancy and distance to park is estimated using the exponential function. But before calibrating such a relationship, we are inspired by Leclercq et al. (2017) to first explore the effects of the time-varying parking conditions on the relationship. Six scenarios are thus compared in Fig. 12c where "increasing/decreasing" suggests that the average parking occupancy increases/decreases during the period from when a vehicle commences cruising on street and to when it successfully parks. Combining the two cases leads to the "both" scenario. We also distinguish between two average parking occupancy values "avg" and "init". To account for the time-varying nature of parking occupancy, the former is calculated as the average of the two average parking occupancy values obtained when a vehicle commences and finishes the on-street parking search. The latter simply refers to the value obtained at the beginning of the search. The cause for the irregular tails of the curves within the rectangles is the very few simulated data lying in the corresponding occupancy ranges[22].

---

[22] For example, the downward trend at the high occupancy end arises from the fact that many vehicles are still cruising at the end of the simulation whose distances to park are not considered.

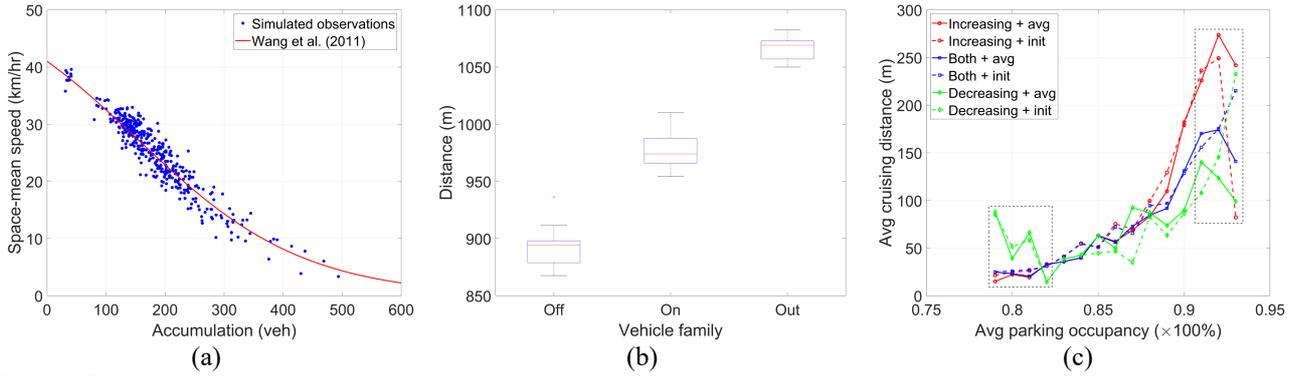

**Fig. 12.** Calibrating the macroscopic model using results obtained from the microscopic model: (a) speed-accumulation NFD; (b) family-specific average moving distances; and (c) relationship between parking occupancy and distance to park.

When the average parking occupancy increases leading to more competitive parking conditions (i.e., vehicles park at a faster rate than that of re-departure), the average distance to park increases accordingly as fewer and fewer spots are available during the search. In contrast, when the average parking occupancy decreases resulting from faster vehicles' re-departures, more and more spots are available during the search leading to a shorter distance to park. Results suggest that, for the same average parking occupancy, the average distance to park varies depending on how the occupancy is evolving. When the parking and re-departure rates are similar, on-street parking is at a quasi-equilibrium state where the average distance to park can be characterized by the blue curve. However, if the parking rate is generally higher/lower than that of re-departure, the average parking occupancy would keep increasing/decreasing requiring that the average distance to park be characterized by the red/green curve. Thus, in a general setting with both peak and off-peak demand, one shall identify the points in time when the trend in the average parking occupancy reverses.

The relationship does not change significantly as we switch between "avg" and "init". Despite slight variations, the solid and dashed curves share a common trend irrespective of the average parking occupancy being "increasing", "decreasing", or "both". The importance of "init" lies in the fact that in the macroscopic model, the average distance to park for any points in time can only be estimated using the current average parking occupancy (i.e. "init"). It is impossible to know a priori the future parking occupancy to calculate "avg". Thus, the calibrated relationship under the scenario "increasing + init" reads $l_c^{on}(O_{on}(t_k)) = 5.2 \times 10^{-11} \exp(24.4 O_{on}(t_k))$ km. We emphasize that using "init" may no longer be valid in the presence of fast-varying parking conditions. Under such circumstances, the average distance to park can be much longer/shorter than that calculated using "init", because the average parking occupancy may have changed significantly during the search. How to characterize distance to park under fast-varying parking conditions requires further research.

*4.2.2. Validation*

Although the macroscopic model is deterministic and parsimonious in nature, it produces encouraging results after calibration that are consistent with those obtained from the stochastic and fine-grained microscopic model. We observe that the accumulation of vehicles in the off-street parking lot keeps increasing until reaching the capacity of 100 at about 2,000 seconds (see Fig. 13a). It remains close to capacity thereafter to the end of the simulation. The accumulation of vehicles parked on street also keeps increasing, but it exhibits a slightly downward trend at the end of the simulation[23] (see Fig. 13b). Both the accumulation of active vehicles in the network and the resulting network speed (see Fig.

---

[23] The macroscopic model slightly overestimates at the peak of the curve, but the relative error is only less than five percent given that the magnitude of the accumulation is over 1,000.

13c and d) can be accurately characterized by the macroscopic model, which lie amid the fluctuations of the stochastic microsimulation results. In some sense, the results of the macroscopic model capture the mean trend of those obtained from the microscopic model.

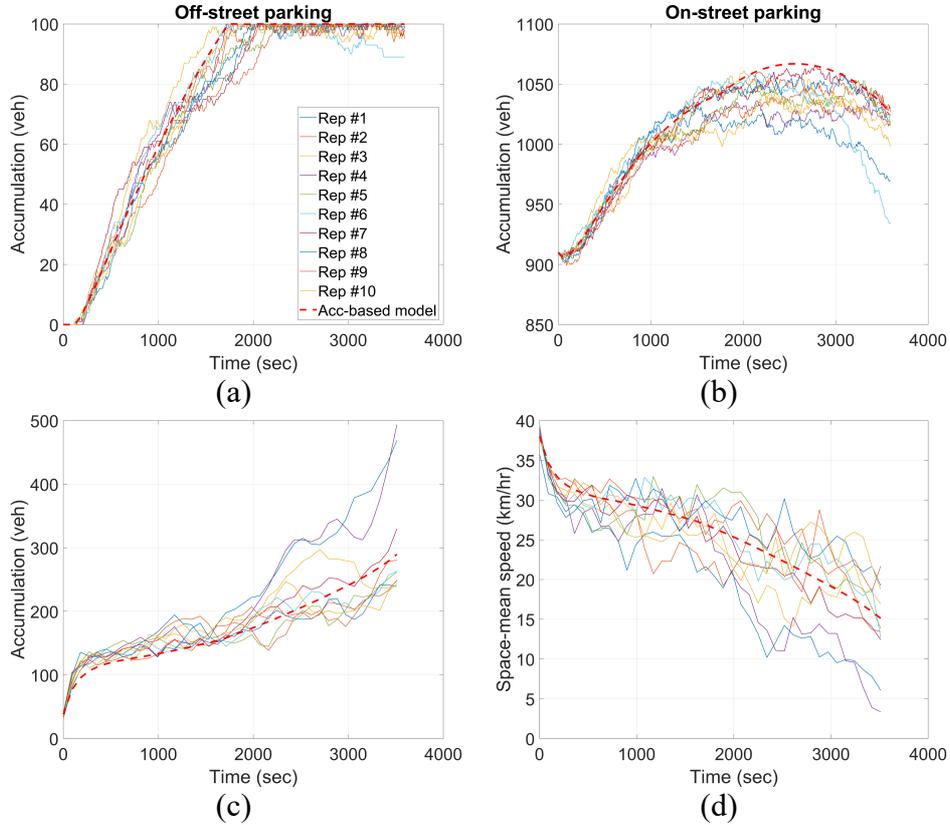

**Fig. 13.** Validating the calibrated macroscopic model: (a) off-street parking accumulation; (b) on-street parking accumulation; (c) accumulation of active vehicles in the network; and (d) space-mean speed of the network.

*4.3. A macro-micro approach to online parking pricing optimization*

Despite its parsimony, the calibrated macroscopic model manifests consistency with the fine-grained microscopic model. It is able to macroscopically characterize the integrated road-parking system dynamics that is known to involve many behavioral complexities of the parking search. The consistency between the models permits integration of the two for online parking pricing optimization, which is to be achieved using model predictive control (MPC). We assume both on- and off-street parking are publicly operated by a central agency with no parking pricing competition. Private parking is not accounted for as it requires introducing another party in the model as well as further behavioral and economic complexities. We refer to Inci and Lindsey (2015); Zheng and Geroliminis (2016) for insights into a parking market with multiple parties.

*4.3.1. Model predictive control*

MPC is a rolling horizon approach that iteratively solves an open-loop optimization problem for a finite prediction horizon and applies the first optimized decision variable to the plant (which is either a real-world network or a simulation model) for the upcoming time steps. Since the model intermittently communicates with the plant, the resulting closed-loop structure is able to address discrepancies between the two due to disturbances to or uncertainties of the network traffic. We refer to Geroliminis et al. (2013) for a brief summary of MPC applications in traffic management and control.

The overall MPC framework is illustrated in Fig. 14 integrating the macro- and micro-models for online parking pricing optimization. Since the time step size of the prediction model $\Delta t'$ can be different from that of the plant $\Delta t$[24], we introduce a new $k'$ in addition to $k$. MPC is not implemented at every time step but only when a control time step $K$ is fully covered. The open-loop optimization problem is solved using the up-to-date family-specific accumulations. The size of $K$, defined as $\Delta T$, is equal to the length of individual parking pricing intervals. Given the 1-hr simulation horizon and the findings from Gu et al. (2020), we choose this length to be 15 minutes considering user adaptation and pricing effectiveness. The number of pricing intervals to be optimized simultaneously is set at two based on the results of Geroliminis et al. (2013), resulting in a 30-min control horizon.

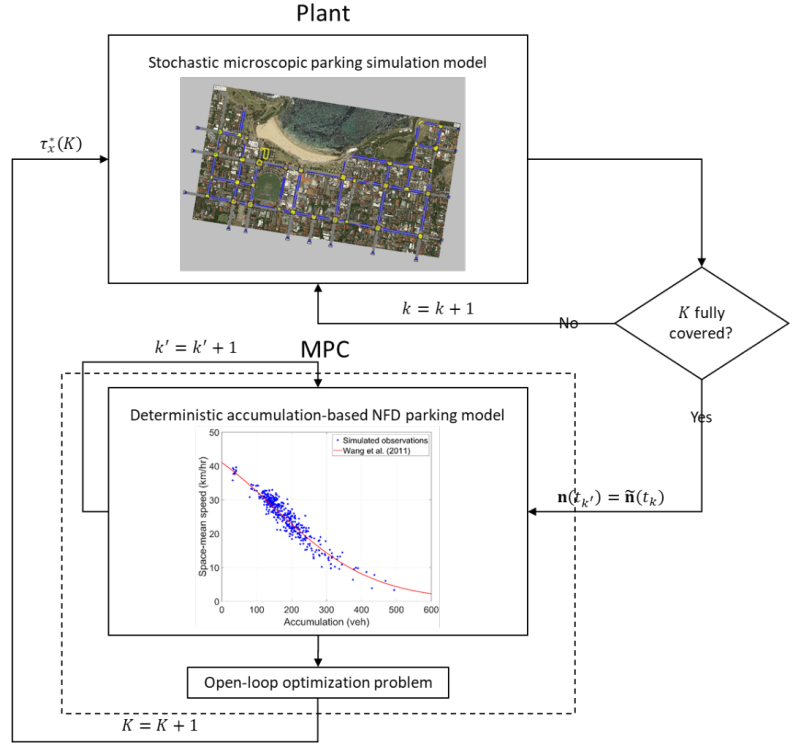

**Fig. 14.** Integrating the macro- and micro-models using the MPC approach for online parking pricing optimization.

The macroscopic model allows using different objectives in the open-loop optimization problem depending on the application needs. Here, we choose to minimize the total ineffective cruising time in order to reduce as much as possible the negative effects of cruising-for-parking on the network traffic. It consists of two quantities. The first quantity is the total cruising time spent by all the vehicles cruising on street $\sum_{k=0}^{\left[\frac{T'}{\Delta t'}\right]-1} n_\mathrm{c}(t_k)\Delta t'$, where $T'$ is the prediction horizon chosen to be 30 minutes[25] comprising $\left[\frac{T'}{\Delta t'}\right]$ time steps. The second quantity is the total cruising time spent by all the vehicles in the off-street parking lot that fail to park $\sum_{k=0}^{\left[\frac{T'}{\Delta t'}\right]-1} \frac{l_\mathrm{off} q_\mathrm{off}^\mathrm{on}(k+1)}{v_\mathrm{off}^\mathrm{f}}$ (which is a pure deadweight loss). We do not include the cruising time of any vehicles that successfully park off street because (i) as long as vehicles manage to park off street, their cruising time is not considered as being ineffective; and (ii) vehicles that successfully park off street do not interact directly with the on-street traffic.

The overall open-loop optimization problem is thus formulated as follows:

---

[24] The microsimulation resolution is 0.2 seconds while the macroscopic model has a time step size of 10 seconds.
[25] Geroliminis et al. (2013) chose the prediction interval to be much longer than the control interval, because the length of their individual control intervals was only 1 min. Here, the length of individual pricing intervals is 15 min. Developing a systematic tuning procedure is not within the scope of the paper.

$$\min_{\boldsymbol{\tau}_x} \sum_{k=0}^{\left[\frac{T'}{\Delta t'}\right]-1} n_c(t_k)\Delta t' + \frac{l_{\text{off}} q_{\text{off}}^{\text{on}}(k+1)}{v_{\text{off}}^{\text{f}}}, \quad (19a)$$

s.t.

Eqs. (10)-(18d),

$$|\tau_x(K+1) - \tau_x(K)| \leq \tau_x^{\text{gap}}, \ K \in \left\{1, \ldots, \left[\frac{T'}{\Delta T}\right]\right\}, \quad (19b)$$

$$\tau_x^{\min} \leq \tau_x(K) \leq \tau_x^{\max}, \ K \in \left\{1, \ldots, \left[\frac{T'}{\Delta T}\right]\right\}, \quad (19c)$$

where $\left[\frac{T'}{\Delta T}\right]$ calculates the number of parking pricing intervals in the prediction horizon, $\boldsymbol{\tau}_x = \left[\tau_x(1), \ldots, \tau_x\left(\left[\frac{T'}{\Delta T}\right]\right)\right]^{\text{T}}$ is the decision vector of the parking price whose first entry $\tau_x(1)$ is to be applied to the plant after optimization, $\tau_x^{\text{gap}}$ is a smoothing parameter chosen to be three dollars in order to prevent unduly fluctuations between adjacent parking prices[26], and $\tau_x^{\min}$ and $\tau_x^{\max}$ are the lower and upper bounds on the parking price chosen to be zero and ten dollars, respectively. To solve the open-loop optimization problem, existing nonlinear optimization methods can be employed (e.g. trust region methods). We further apply a multi-start technique to prevent the optimization from being trapped in a (bad) local optimum given the non-smooth objective function. However, although the solution quality is improved at the cost of more computational time, the global optimum is not strictly guaranteed.

*4.3.2. Scenario analysis: on-street parking pricing*

There are concerns among transportation economists about the adverse effects of free or under-priced on-street parking (Shoup, 2006). To provide computational evidence in this regard (see Section D in the *Supplementary Material* for an alternative scenario analysis on off-street parking pricing), we assume the capacity of the off-street parking lot is doubled so that it becomes underutilized while vehicles are still competing for limited on-street parking resources. Two additional scenarios, "full-horizon dynamic" and "full-horizon static", are further considered in the comparison, where a single full-horizon optimization problem is formulated and solved. The difference is that the price is time-varying in the former (a total of four prices) and not in the latter.

Results show that the optimal price profiles of all the ten replications share a common pattern, which is also observed under the "full-horizon dynamic" scenario (see Fig. 15a). Under the "full-horizon static" scenario, the optimal price profile is a straight line lying amid the other curves, as expected. When price varies with time, the peak occurs at the first interval driving vehicles to park in the off-street parking lot; otherwise, on-street parking would trigger a significant cruising-for-parking effect. As the off-street parking lot becomes increasingly occupied, the price gradually declines until reaching the third interval and then rebounds at the end of the simulation. We illustrate in Fig. 15b the capability of the macroscopic model in predicting and correcting the network state in real-time according to the microsimulation results (also see Section C in the *Supplementary Material*).

In Fig. 15c we compare four time-related network performance metrics under four scenarios consisting of "no-price", "MPC", and the two additional scenarios previously described. The four metrics feature (a) off-street cruising deadweight loss; (b) on-street total cruising time; (c) total ineffective cruising time (i.e. the sum of (a) and (b)); and (d) total travel time on road. When on-street parking is optimally priced, more vehicles intend to use the off-street parking lot resulting in a higher competition and thus a higher probability of failing to park. This explains the slight increase in the off-street deadweight loss (see panel (a) in Fig. 15c). The increase is indeed slight compared with the decrease in the total on-street cruising time thanks to fewer vehicles cruising on street (see panel (b) in Fig. 15c). Thus, we achieve an overall reduction in the total ineffective cruising time (see panel (c) in Fig. 15c). The most counterintuitive observation lies in the total on-road travel time which does not reduce much as expected (see panel (d) in Fig. 15c). The cause is found to be the roundabout located upstream of the off-street parking lot, which fails to provide sufficient capacity to accommodate the increased number of vehicles intending to park off street. Thus, vehicles queue at the roundabout creating a local congestion effect that cancels out the time savings arisen from fewer vehicles cruising on street[27].

---

[26] Such constraints can also be found in Geroliminis et al. (2013) and Gu et al. (2019).

To provide further simulated evidence, we replace the roundabout with a hypothetical interchange so that vehicles coming from each approach have their respective right of way without any conflicts. This is the only difference between the original and the modified networks. In the original network, we have already seen that the total on-road travel time does not reduce much even though on-street parking is optimally priced. This is no longer the case in the modified network where a clear reduction is present (see Fig. 15d). The result emphasizes the importance of accessibility to off-street parking when on-street parking is priced. In other words, the road network leading to off-street parking must have sufficient capacity to accommodate vehicles that are priced away from parking on street.

We emphasize that the above observation would become more notable if we were to minimize both the total on-road travel time and the off-street deadweight loss, because (i) the optimal price would become higher to discourage more vehicles from parking on street[28]; and (ii) more vehicles intending to park off street would lead to more vehicles failing to park that further congest the roundabout. Given the principle of proximity, adverse effects could ensue including queue spillback and partial network gridlock. While this worst situation does not arise in the simulation, in reality we do see vehicles queue and spill back onto streets nearby off-street parking. From a theoretical perspective, the violation of the homogeneity assumption underlying the macroscopic approach requires that the network be partitioned into multiple small parking areas. This is of future research priority.

---

[27] Pricing the congested city center without providing sufficient capacity to accommodate the detour traffic would similarly result in queues along the periphery contributing to local congestion (Gu et al., 2018b).
[28] Fewer vehicles cruising on street help reduce the travel time of those in transit.

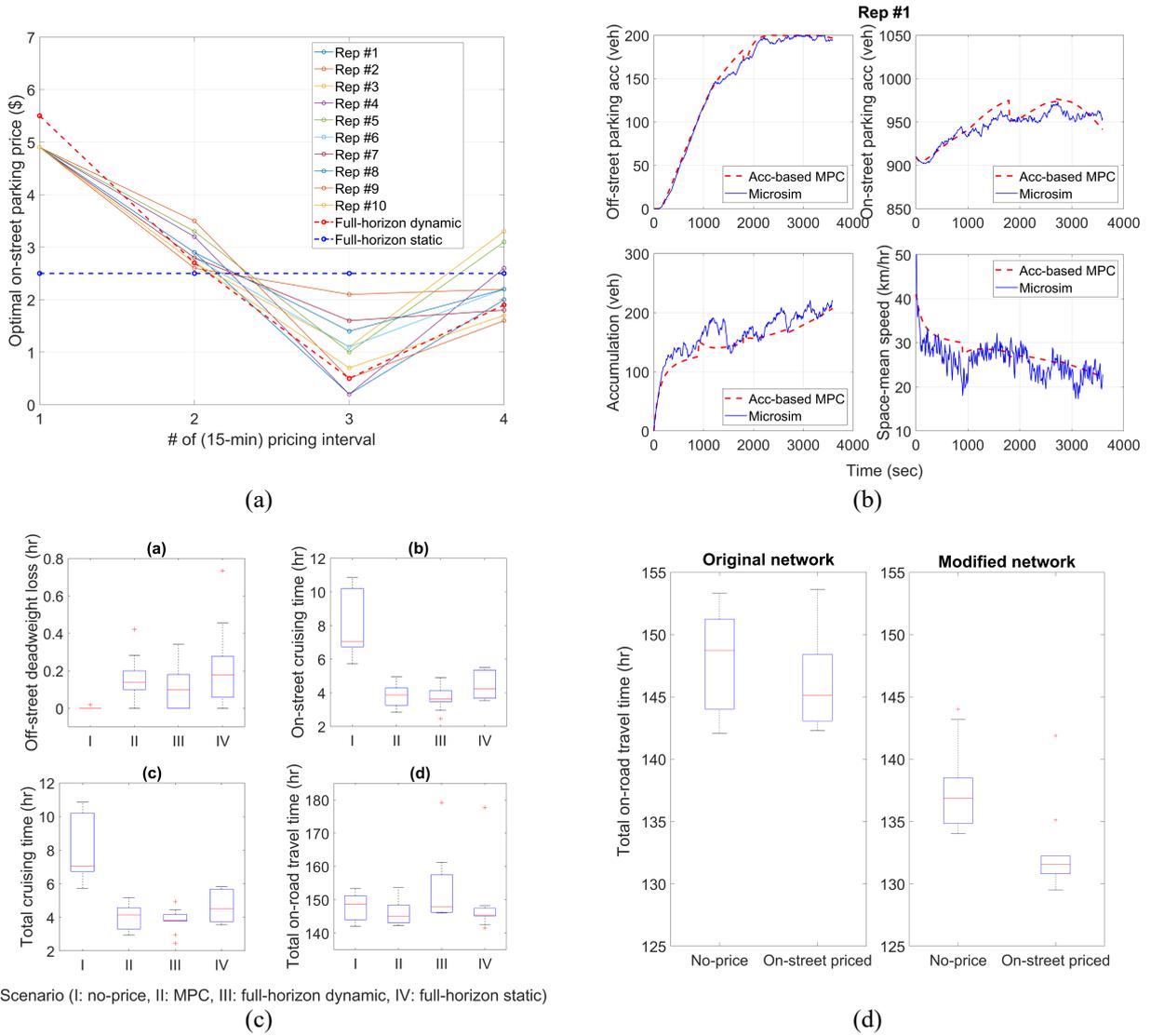

**Fig. 15.** Applying the MPC approach to optimize the on-street parking price in real-time: (a) optimal price profiles; (b) comparative results of the macro- and micro-models for replication 1; (c) network performance metrics comparison; and (d) total on-road travel time comparison between the original and the modified networks.

## 5. Conclusion

### 5.1. Summary of findings

Understanding and modeling parking dynamics is known to be challenging due to many behavioral complexities such as cruising-for-parking. In this paper, we develop a microscopic parking simulation model that explicitly considers cruising-for-parking, based on which some key aspects of parking modeling are discussed:

i. Although the low cruising speed reduces the network performance, it does not significantly alter the NFD unless cruising vehicles dominate the traffic stream. Thus, the simulated or empirical NFD can usually be used without the need to account for the effects of low cruising speed.

ii. Distance to park is not uniquely determined by parking occupancy because factors such as cruising speed and parking duration also contribute. In general, a lower cruising speed results in a shorter distance to park while a longer parking duration results in a longer distance to park.
  iii. Multiscale parking occupancy-driven intelligent parking guidance can effectively reduce the spatial heterogeneity of the parking pattern, resulting in shorter distances to park and thus considerable network efficiency gains.

Using the microscopic model, we extend, calibrate, and validate a macroscopic parking dynamics model with an NFD representation. Despite being parsimonious, the calibrated macroscopic model is able to capture the integrated road-parking system dynamics arisen from the much more detailed microscopic model, which permits integration of the two resulting in a new macro-micro approach to modeling parking. We apply the proposed approach to online parking pricing optimization, and the results demonstrate its effectiveness. However, one caveat of the approach is that, when pricing on-street parking, the road network connected to the alternate off-street parking lots must have sufficient capacity to accommodate the increased parking demand; otherwise, local congestion may arise that violates the homogeneity assumption underlying the macroscopic model.

*5.2. Future research*

Although this paper provides insights into multiple aspects of parking modeling, many other questions remain open and thus require further research efforts. How to accurately estimate distance to park is one such question. Multiple factors contribute including, among others, parking occupancy, parking duration, cruising speed, and parking competition. Thus, a systematic study is encouraged to reveal the relationship between these factors and distance to park. A natural extension to the current work is to propose, calibrate, and validate a multi-reservoir parking dynamics model. How to partition a network into multiple parking areas is one question to be answered, as existing studies on network partitioning (Gu and Saberi, 2019; Ji and Geroliminis, 2012; Saeedmanesh and Geroliminis, 2016) did not consider parking in their networks. The other two questions are how to model the parking area choices of individuals (Boyles et al., 2015) and the resulting flow transfers between neighboring areas (Mariotte and Leclercq, 2019). Finally, the alternate trip-based approach to modeling parking is of future research priority given its capability of incorporating user heterogeneity such as individualized trip length (Leclercq et al., 2017).

**Acknowledgments**

The authors acknowledge the financial support from the Department of Industry, Innovation and Science under the Grant Agreement SCS69276. Rashidi also acknowledges the financial support from the Australian Research Council (DE170101346). The authors appreciate the constructive suggestions from three anonymous referees that significantly improve the quality of the paper.